 \documentclass[10pt]{article}

\usepackage[backend=biber, style=alphabetic, citestyle=alphabetic, backref=true]{biblatex}
\usepackage[margin=1in]{geometry}
\usepackage[usenames,dvipsnames,svgnames,table]{xcolor}
\usepackage{graphicx, float, wrapfig, caption, subcaption}

\usepackage{fancyhdr}
\usepackage{enumerate}

\usepackage{yfonts}
\usepackage[mathscr]{eucal}

\usepackage{etoolbox}

\newbool{ShowFootnotes}

\newbool{ShowFutureWork}

\usepackage[final]{hyperref}
\hypersetup{colorlinks=true}
\hypersetup{linkcolor=Emerald}
\hypersetup{citecolor=cyan}
\hypersetup{menucolor=red}

\usepackage{beesonmacros}
\newcommand{\eps}{\epsilon}
\newcommand{\Xe}{X^\epsilon}
\newcommand{\Ze}{Z^\epsilon}
\newcommand{\Ye}{Y^\epsilon}
\newcommand{\Me}{M^\epsilon}

\newcommand{\mYe}{\mathcal{Y}^\epsilon}
\newcommand{\mF}{\mathcal{F}}
\newcommand{\wt}[1]{\widetilde{#1}}
\newcommand{\Pe}{\mathbb{P}^\epsilon}
\newcommand{\mG}{\mathcal{G}}

\newcommand{\wDe}{\wt{D}^\epsilon}

\newcommand{\later}[1]{\textcolor{Fuchsia}{#1}}
\newcommand{\FutureWork}[1]{\ifbool{ShowFutureWork}{\later{#1}}{}}
\newcommand{\Nicolas}[1]{\textcolor{black}{#1}}

\newcommand{\orcid}[1]{\href{https://orcid.org/#1}{\includegraphics[scale=.014]{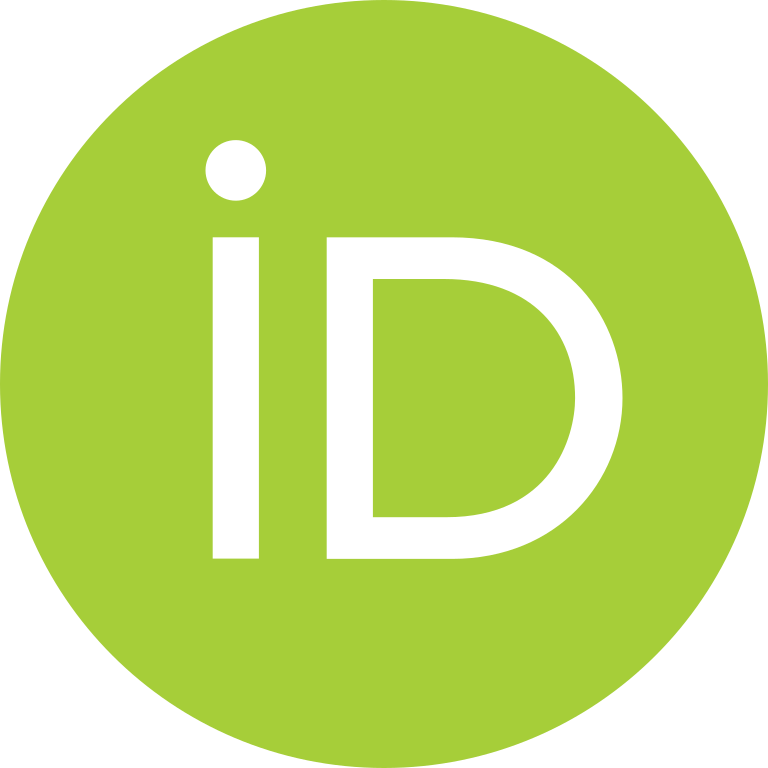}}}

\title{Approximation of the Filter Equation for Multiple Timescale, Correlated, Nonlinear Systems}
\date{\today}
\author{Ryne Beeson\footnote{University of Illinois at Urbana-Champaign} \orcid{0000-0003-2176-0976},
N. Sri Namachchivaya\footnote{University of Waterloo}, 
and Nicolas Perkowski\footnote{Freie Universit\"at Berlin}}

\begin{document}

\maketitle

\begin{abstract}
This paper considers the approximation of the continuous time filtering equation for the case of a multiple timescale (slow-intermediate, and fast scales) that may have correlation between the slow-intermediate process and the observation process. 
The signal process is considered fully coupled, taking values in $\R{m} \times \R{n}$ and without periodicity assumptions on coefficients. 
It is proved that in the weak topology, the solution of the filtering equation converges in probability to a solution of a lower dimensional averaged filtering equation in the limit of large timescale separation.
The method of proof uses the perturbed test function approach (method of corrector) to handle the intermediate timescale in showing tightness and characterization of limits. 
The correctors are solutions of Poisson equations. 
\end{abstract}


\section{Introduction}

The aim of this paper is to prove a convergence result for the continuous time filtering equation of a multiple timescale and correlated nonlinear system to a lower dimensional filtering equation. 
Specifically, consider the coupled system of stochastic differential equations (SDEs), 
\EquationAligned{
\label{equation: multiscale signal process}
d\Xe_t 
&= \left[ b(\Xe_t, \Ze_t) + \frac{1}{\epsilon} b_I(\Xe_t, \Ze_t) \right] dt + \sigma(\Xe_t, \Ze_t) dW_t, \\
d\Ze_t 
&= \frac{1}{\eps^2} f(\Xe_t, \Ze_t) dt + \frac{1}{\eps}g(\Xe_t, \Ze_t) dV_t,
}
and denote the infinitesimal generator of $(\Xe, \Ze)$ as $\mG^\epsilon$.
$(\Xe, \Ze)$ is known as the signal process and $\epsilon \in (0, 1)$ is a timescale parameter such that $\Ze$ is a fast process and $\Xe$ is a slow process. 
Note that even the equation for $\Xe$ possesses an intermediate timescale due to the $\frac{1}{\epsilon} b_\mathrm{I}$ drift coefficient. 
In filtering theory, we consider the signal process to be non-observable, and instead have indirect measurements of $(\Xe, \Ze)$ via the noisy observation process
\begin{align*}
d\Ye_t = h(\Xe_t, \Ze_t) dt + \alpha dW_t + \gamma dU_t.
\end{align*}
With $W, V, U$ independent Brownian motions, $\alpha \neq 0$ indicates correlation between the slow signal $\Xe$ and the observation process $\Ye$.  
The goal in filtering theory is then to calculate the conditional distribution of $(\Xe, \Ze)$ given the observation history generated from $\Ye$, which we denote by $\pi^\epsilon$. 
At each time $t > 0$, $\pi^\epsilon_t$ is a random probability measure on the space $\R{m} \times \R{n}$ and acts on test functions $\varphi : \R{m} \times \R{n} \rightarrow \R{}$ by integration $\pi^\epsilon_t(\varphi) = \int \varphi(x, z) \pi^\epsilon_t(dx, dz)$. 

The motivating question of this paper then comes from the known result that if for every fixed $x$, the solution $Z^x$ of
\begin{align*}
dZ^x_t = f(x, Z^x_t) dt + g(x, Z^x_t) dV_t,
\end{align*}
is ergodic with stationary distribution $\mu_\infty(x)$, then under appropriate assumptions, the process $\Xe$ converges in distribution to a Markov process $X^0$ with infinitesimal generator $\mG^\dagger$ in the limit as $\epsilon \rightarrow 0$ \cite{Papanicolaou:1977dtc, Pardoux:2003, Khasminskii:2005jde}.
Therefore, if we are only interested in statistics of $\Xe$ (i.e. estimation of test functions $\varphi : \R{m} \rightarrow \R{}$), then it would be computationally advantageous to know if $\pi^{\epsilon, x} \Rightarrow \pi^0$ converges weakly to a lower dimensional filtering equation; $\pi^0_t$ being a random probability measure for each time $t$ on $\R{m}$ and $\pi^{\epsilon, x}$ being the $x$-marginal of $\pi^\epsilon$. 

Filtering theory has widespread applications in many fields including various disciplines of engineering for decision and control systems, the geosciences, weather and climate prediction. 
In many of these fields, it is not uncommon to have physics based models with multiple timescales as seen in Eq. \ref{equation: multiscale signal process}, and also have the case were estimation of the slow process is solely of interest; for example the estimation of the ocean temperature, which is necessary for climate prediction, but the ocean model may also be coupled to a fast atmospheric model. 
Knowing that mathematically $\pi^{\epsilon, x} \Rightarrow \pi^0$ in the limit as $\epsilon \rightarrow 0$, enables practitioners to devise more efficient methods for estimation of the slow process without great loss of accuracy (see for instance \cite{Park:2011jam, Kang:2012hj, Berry:2014, Yeong:2020}).

There are several papers providing results for $\pi^{\epsilon, x} \rightarrow \pi^0$ (or the associated unnormalized conditional measure or density versions) on variations of the multiple timescale filtering problem. 
In \cite{Park:2010n}, $(\Xe, \Ze)$ is a two dimensional process with no drift in the fast component, no intermediate scale, and no correlation. 
The authors made use of a representation of the slow component by a time-changed Brownian motion under a suitable measure to yield weak convergence of the filter. 
Homogenization of the nonlinear filter was studied in \cite{Bensoussan:1986s} and \cite{Ichihara:2004} by way of asymptotic analysis on a dual representation of the nonlinear filtering equation. 
In these papers, the coefficients of the signal processes are assumed to be periodic. 
The approach in \cite{Ichihara:2004} is novel as the first application of backward stochastic differential equations for homogenization of Zakai-type stochastic partial differential equations (SPDEs). 

Convergence of the filter for a random ordinary differential equation with intermediate timescale and perturbed by a fast Markov process was investigated in \cite{Lucic:2003amo}. 
A two timescale problem with correlation between the slow process and observation process, but where the slow dispersion coefficient does not depend on the fast process, 
is investigated in \cite{Qiao:2019o}.
The main result is that the filter converges in $L^1$ sense to the lower dimensional filter. 
An energy method approach is used in \cite{Zhang:2019sd} to show that the probability density of the reduced nonlinear filtering problem approximates the original problem when the signal process has constant diffusion coefficients, periodic drift coefficients and the observation process is only dependent on the slow process. 

Convergence of the nonlinear filter is shown in a very general setting in \cite{Kleptsina:1997}, based on convergence in total variation distance of the law of $(\Xe, \Ye)$. 
In the examples of \cite{Kleptsina:1997}, the diffusion coefficient is not allowed to depend on the fast component.

In contrast to other papers on the convergence of the nonlinear filter for the multiple timescale problem, Imkeller et al. \cite{Imkeller:2013} showed a quantitative rate of convergence of $\epsilon$ for the system in Eq. \ref{equation: multiscale signal process}, but without intermediate timescale nor correlation of the slow process with the observation process. 
This is accomplished using a suitable asymptotic expansion of the dual of the Zakai equation and then harnessing a probabilistic representation of the SPDEs in terms of backward doubly stochastic differential equations. 
This result is then extend to the case of correlation between the slow signal and observation process in \cite{Beeson:2020jco}, with the same rate of convergence. 

We lastly mention the work of nonlinear filter approximation given in \cite[Chapter 6]{Kushner:1990}, which is most similar to the approach used in  this paper. 
In \cite[Chapter 6]{Kushner:1990}, a two timescale jump-diffusion process is considered, but with no correlation between signal and observation process. 
The difference of the actual unnormalized conditional measure and the reduced conditional measure converges to zero in distribution.
Standard results then yield convergence in probability of the fixed time marginals. 
The method of proof is by averaging the coefficients of the SDEs for the unnormalized filters and showing that the limits of both filters satisfy the same SDE, which possess a unique solution.  

In this paper, we address the broader multiple timescale correlation filtering problem and therefore have to modify the approach by \cite[Chapter 6]{Kushner:1990} to handle the intermediate scaling term and the correlation. 
This is the first paper that the authors are aware of that handles the problem of filter convergence for correlated slow-fast systems with intermediate timescale forcing. 
For this we make use of the perturbed test function approach where the correctors are solutions of Poisson equations. 
We make use of the sharp results on existence, regularity and growth of the transition densities and semigroups associated with the process $Z^x$ and Poisson equations for the corrector terms in \cite{Pardoux:2003}.
The main result of the paper is the following:
\begin{theorem*}[Main Result]
Recall that $\pi^{\epsilon, x}$ is the $x$-marginal of the conditional distribution $\pi^\epsilon$ and $\pi^0$ is the conditional distribution for the averaged filter equation (see for instance Eqs. \ref{equation: normalized conditional distribution as expectation}, \ref{equation: marginal conditional distribution as expectation}, and \ref{equation: homogenized kushner-stratonovich}). 
Then under the assumptions stated in Theorem \ref{theorem: main result, weak convergence of the normalized filter}, 
$\pi^{\epsilon, x} \rightarrow \pi^0$ in probability.
\end{theorem*}

The paper proceeds as follows: 
Section \ref{section: problem statement} provides the problem statement in greater detail and states the main theorem with full assumptions.
We also introduce the unnormalized variants of $\pi^\epsilon$ and $\pi^0$, which we will denote as $\rho^\epsilon$ and $\rho^0$ in this section.
Similar to the notation $\pi^{\epsilon, x}$, $\rho^{\epsilon, x}$ will denote the $x$-marginal of $\rho^\epsilon$. 
In Section \ref{section: preliminary estimates} we give preliminary estimates which are needed for the main results. 
Section \ref{section: existence, characterization, and uniqueness} provides the existence of weak limits of the probability measure induced by the signed measured-valued process $\rho^{\epsilon, x} - \rho^0$, as well as the characterization of this process and the uniqueness of its limit. 
At the end of Section \ref{section: existence, characterization, and uniqueness}, the main result for convergence of $\rho^{\epsilon, x} - \rho^0$ is stated alongside a lemma that proves the convergence of $\pi^{\epsilon, x} - \pi^0$.

\section{Problem Statement}
\label{section: problem statement}

In this section, we provide the full problem statement, some notation and the main result. 
We consider a filtered probability space $(\Omega, \mathcal{F}, (\mathcal{F}_t)_{t \geq 0}, \Q{})$ supporting a $(w + v + u)$-dimensional $\mF_t$-adapted Brownian motion $(W, V, U)$.
We will work with the following system of SDEs,
\EquationAligned{
\label{equation: original SDE problem setup}
dX^\epsilon_t 
&= \left[ b(X^\epsilon_t, Z^\epsilon_t) + \frac{1}{\epsilon} b_I(X^\eps_t, Z^\eps_t) \right] dt + \sigma(X^\epsilon_t, Z^\epsilon_t) dW_t, \\[1.5mm]
dZ^\eps_t 
&= \frac{1}{\eps^2} f(X^\eps_t, Z^\eps_t) dt + \frac{1}{\eps}g(X^\eps_t, Z^\eps_t) dV_t, \\[3mm]
dY^\eps_t 
&= 
h(X^\eps_t, Z^\eps_t) dt + \alpha dW_t + \gamma dU_t,  \quad \Ye_0 = 0 \in \R{d},\\
}
where $b, b_\mathrm{I} : \R{m} \times \R{n} \rightarrow \R{m}$, $\sigma : \R{m} \times \R{n} \rightarrow \R{m} \times \R{w}$, $f : \R{m} \times \R{n} \rightarrow \R{n}$, $g : \R{m} \times \R{n} \rightarrow \R{n} \times \R{v}$ and $h : \R{m} \times \R{n} \rightarrow \R{d}$ are Borel measurable functions. 
The initial distribution of $(X, Z)$ is denoted by $\Q{}_{(\Xe_0, \Ze_0)}$ and is assumed independent of the $(W, V, U)$ Brownian motion.
$\Q{}_{(\Xe_0, \Ze_0)}$ is also assumed to have finite moments for all orders\ifbool{ShowFootnotes}{\footnote{appears that we could relax to just finite moments of up to and including second order as needed in (i) of Lemma \ref{lemma: solutions to signed Zakai equations are tight}}}{}. 
In Eq. \ref{equation: original SDE problem setup}, $0 < \epsilon \ll 1$, is a timescale separation parameter.
We consider the case where $\alpha \in \R{d \by w}, \gamma \in \R{d \by u}$, and assume the following to be true
\begin{align*}
K \equiv \alpha \alpha^* + \gamma \gamma^* \succ 0, \quad \gamma \gamma^* \succ 0.
\end{align*}
This implies the existence of a unique $\R{d \by d} \ni \kappa \succ 0$ of lower triangular form, such that $K = \kappa \kappa^*$. 
Hence there exists a unique $\kappa^{-1}$, such that we can define an auxiliary observation process
\begin{align}
\label{equation: auxiliary observation process}
Y^{\epsilon, \kappa}_t = \int_0^t \kappa^{-1} d\Ye_s = \int_0^t \kappa^{-1} h(\Xe_s, \Ze_s) ds + B_t, \quad Y^{\epsilon, \kappa}_0 = 0 \in \R{d},
\end{align}
where
\begin{align*}
B_t = \kappa^{-1} \left( \alpha dW_t + \gamma dU_t \right), 
\end{align*}
is a standard $d$-dimensional Brownian motion under $\Q{}$. 

We are interested in the convergence of the normalized filter, $\pi^\epsilon$, the conditional distribution of the signal given the observation filtration, to an averaged form. 
In particular, for any test function $\varphi \in C^2_b(\R{m} \times \R{n}; \R{})$ and time $t \in [0 , T]$, the normalized filter can be characterized as
\begin{align}
\label{equation: normalized conditional distribution as expectation}
\pi^\epsilon_t(\varphi) = \E[\Q{}]{ \varphi(\Xe_t, \Ze_t) }[\mYe_t], 
\end{align}
where $\mYe_t \equiv \sigma(\set{\Ye_s}[s \in [0, t]]) \vee \mathcal{N}$, the $\sigma$-algebra generated by the observation process over the interval $[0, t]$, joined with $\mathcal{N}$, the $\Q{}$ negligible sets. 

Because the filtrations generated by $\Ye$ and $Y^{\epsilon, \kappa}$ are equivalent, from the point of view of $\pi^\epsilon$ we can use either. 
Hence, let us redefine the sensor function $h \leftarrow \kappa^{-1} h$, and $\alpha \leftarrow \kappa^{-1} \alpha$, $\gamma \leftarrow \kappa^{-1} \gamma$, so that the observation process can be redefined as
\begin{align}
\label{equation: redefined observation process}
dY^\eps_t = h(X^\eps_t, Z^\eps_t) dt + dB_t,  \quad \Ye_0 = 0 \in \R{d},
\end{align}
where $B = \alpha W + \gamma U$ is a standard Brownian motion under $\Q{}$ and still correlated with $W$. 

In Eq. \ref{equation: original SDE problem setup}, we identify the infinitesimal generators of the SDEs as follows, 
\begin{align*}
\mG_S (x, z) &\equiv \sum_{i=1}^m b_i(x, z) \frac{\partial}{\partial x_i} + \frac{1}{2} \sum_{i, j=1}^m (\sigma \sigma^*)_{ij} (x, z) \frac{\partial^2}{\partial x_i \partial x_j}, \\
\mG_I (x, z) &\equiv \sum_{i=1}^m b_{I, i} (x, z) \frac{\partial}{\partial x_i}, \\
\mG_F (x, z) &\equiv \sum_{i=1}^n f_i(x, z) \frac{\partial}{\partial z_i} + \frac{1}{2} \sum_{i, j=1}^n (gg^*)_{ij}(x, z) \frac{\partial^2}{\partial z_i \partial z_j},  \\
\mG^\epsilon_S &\equiv \frac{1}{\epsilon} \mG_I + \mG_S, \\
\mG^\epsilon &\equiv \frac{1}{\epsilon^2} \mG_F + \frac{1}{\epsilon} \mG_I + \mG_S.
\end{align*}
The Kushner-Stratonovich equation for the time evolution of the filter $\pi^\epsilon$, acting on a test function $\varphi \in C^2_b(\R{m} \by \R{n}; \R{})$, is
\EquationAligned{\label{equation: kushner-stratonovich}
\pi^\epsilon_t(\varphi) &= \pi^\epsilon_0(\varphi) + \int_0^t \pi^\epsilon_s (\mG^\epsilon \varphi) ds + \int_0^t \langle \pi^\epsilon_s ( \varphi h + \alpha \sigma^* \nabla_x \varphi ) - \pi^\epsilon_s (\varphi) \pi^\epsilon_s (h), d\Ye_s - \pi^\epsilon_s (h) ds \rangle, \\
\pi^\epsilon_0(\varphi) &= \E[\Q{}]{ \varphi(\Xe_0, \Ze_0) }.
}
When we are interested in estimating test functions of $\Xe$ only, i.e. $\varphi \in C^2_b(\R{m}; \R{})$, we consider the $x$-marginal of $\pi^\epsilon$, 
\begin{align}
\label{equation: marginal conditional distribution as expectation}
\pi^{\epsilon, x}_t(\varphi) = \int \varphi(x) \pi^\epsilon_t(dx, dz).
\end{align}

\subsection{Homogenization of Stochastic Differential Equations}

The theory of homogenization of SDEs shows that if the processes $Z^{\epsilon, x}$, 
\begin{align}
\label{equation: SDE, fast process with fixed slow state}
dZ^{\eps, x}_t = \frac{1}{\eps^2} f(x, Z^{\eps, x}_t) dt + \frac{1}{\eps}g(x, Z^{\eps, x}_t) dV_t,
\end{align}
 is ergodic with stationary distribution $\mu_\infty(x)$, then under appropriate conditions, in the limit $\epsilon \rightarrow 0$ the process $\Xe$ converges in distribution to a Markov process $X^0$ with infinitesimal generator, 
\begin{align*}
\mG^\dagger &\equiv \overbar{\mG_S} + \wt{\mG},
\end{align*}
where
\begin{align*}
\overbar{\mG_S} (x) &\equiv \sum_{i = 1}^m \overbar{b}_i(x) \frac{\partial}{\partial x_i} + \frac{1}{2} \sum_{i, j = 1}^m \overbar{a}_{ij}(x) \frac{\partial^2}{\partial x_i \partial x_j}, \\
\overbar{b}(x) &\equiv \int_{\R{n}} b(x, z) \mu_\infty(dz; x), \\
\overbar{a}(x) &\equiv \int_{\R{n}} a(x, z) \mu_\infty(dz; x), 
\end{align*}
$a = \sigma \sigma^*$, and
\begin{align*}
\wt{\mG} (x) &\equiv \sum_{i = 1}^m \wt{b}_i(x) \frac{\partial}{\partial x_i} + \frac{1}{2} \sum_{i, j = 1}^m \wt{a}_{ij}(x) \frac{\partial^2}{\partial x_i \partial x_j}, \\
\wt{b}(x) &\equiv \int_{\R{n}} \left( \nabla_x \mG_F^{-1} (- b_{\mathrm{I}}) \right) b_{\mathrm{I}} (x, z) \mu_\infty(dz; x), \\
\wt{a}(x) &\equiv \int_{\R{n}} \left( b_{\mathrm{I}} \otimes \mG_F^{-1}( - b_{\mathrm{I}}) \right)(x, z) + \left( \mG_F^{-1}( - b_{\mathrm{I}}) \otimes b_{\mathrm{I}}  \right)(x, z) \mu_\infty(dz; x),
\end{align*}
where $\mG_F^{-1} (- b_{\mathrm{I}})$ is the solution of a Poisson equation. 

We define the averaged filter $\pi^0$, a probability measure-valued process satisfying the following evolution equation,
\EquationAligned{
\label{equation: homogenized kushner-stratonovich}
\pi^0_t(\varphi) &= \pi^0_0(\varphi) + \int_0^t \pi^0_s (\mG^\dagger \varphi) ds + \int_0^t \langle \pi^0_s ( \varphi \overbar{h} + \alpha \overbar{\sigma}^* \nabla_x \varphi ) - \pi^0_s (\varphi) \pi^0_s (\overbar{h}),  d\Ye_s - \pi^0_s (\overbar{h}) ds \rangle, \\
\pi^0_0(\varphi) &= \E[\Q{}]{ \varphi(X^0_0) }.
}
The definitions of $\overbar{h}, \overbar{\sigma}$ are 
\begin{align*}
\overbar{h}(x) \equiv \int_{\R{n}} h(x, z) \mu_\infty(dz; x), \quad \overbar{\sigma}(x) \equiv \int_{\R{n}} \sigma(x, z) \mu_\infty(dz; x).
\end{align*}

\subsection{Notation and Main Theorem}

Before stating the main result of the paper, we set and provide a few definitions and assumptions that will be used throughout the paper. 
We will use $\N{}_0$ to denote $\{ 0, 1, 2, \hdots \}$ and $\N{}$ for $\set{1, 2, \hdots}$.
Let $H_f$ denote the assumption that there exists a constant $C > 0$, exponent $\alpha > 0$ and an $R > 0$ such that for all $|z| > R$,  
\begin{align}
\tag{$H_f$}
\label{assumption: positive recurrence}
\sup_{x \in \R{m}} \langle f(x, z), z \rangle \leq - C |z|^\alpha.
\end{align}
\ref{assumption: positive recurrence} is a recurrence condition, which provides the existence of a stationary distribution, $\mu_\infty(x)$, for the process $Z^x$.
Let $H_g$ denote the assumption that there are $0 < \lambda \leq  \Lambda < \infty$, such that for any $(x, z) \in \R{m} \times \R{n}$,
\begin{align}
\tag{$H_g$}
\label{assumption: uniform ellipticity}
\lambda I \preceq g g^*(x, z) \preceq \Lambda I, 
\end{align}
where $\preceq$ is the order relation in the sense of positive semi-definite matrices. 
\ref{assumption: uniform ellipticity} is a uniform ellipticity condition, which provides the uniqueness of the stationary distribution. 
We will say that a function $\theta: \R{m} \times \R{n} \rightarrow \R{}$ is centered with respect to $\mu_\infty(x)$, if for each $x$ 
\[
\int \theta(x, z) \mu_\infty(dz; x) = 0, \quad \forall x \in \R{m}. 
\]
For brevity in the results to follow, let us denote $H^{i, j + \alpha}$ to specify the regularity and boundedness of $f$ and $gg^*$ as follows, 
\hypertarget{H}{}
\begin{align}
\tag{$H^{i, j + \alpha}$}
\label{assumption: regularity of fast coefficients}
f \in C^{i, j + \alpha}_b(\R{m} \times \R{n}; \R{n}), \quad g g^* \in C^{i, j + \alpha}_b(\R{m} \times \R{n}; \R{n \times n}), \quad i, j \in \N{}, \quad \alpha \in (0, 1),
\end{align}
where $\varphi(x, z) \in C^{i, j + \alpha}_b$ denotes that $\varphi$ has $i$ bounded derivatives in the $x$-component, $j$ bounded derivatives in the $z$-component, and all derivatives $\partial^{j'}_z \partial^{i'}_x \varphi$ for $0 \leq i' \leq i$, $0 \leq j' \leq j$ are $\alpha$-H\"older continuous in $z$ uniformly in $x$. 

We use the notation $k = (k_1, \hdots, k_m) \in \mathbb{N}^m_0$ for a multiindex with order $|k| = k_1 + \hdots + k_m$ and define the differential operator
\begin{align*}
D^k_x = \frac{\partial^{|k|}}{\partial {x_1}^{k_1} \hdots \partial x_m^{k_m}}.
\end{align*}

{
\thm{
\label{theorem: main result, weak convergence of the normalized filter}
Assume that $f$ and $g$ satisfy \ref{assumption: positive recurrence} and \ref{assumption: uniform ellipticity}, that $b_\mathrm{I}$ is centered with respect to $\mu_\infty(x)$ for each $x$ and that $\Q{}_{(\Xe_0, \Ze_0)}$ has finite moments of every order.
Additionally, assume either a. high regularity conditions or b. low regularity with uniform ellipticity:
\begin{enumerate}[a.]
\item \hyperlink{H}{$H^{4, 2 + \alpha}$} holds for $\alpha \in (0, 1)$;
for each $z$, $b(\cdot, z) , \sigma(\cdot, z) \in C^3$, and $b_\mathrm{I}(\cdot, z) \in C^4$;
that $b$ and $b_\mathrm{I}$ are Lipschitz in $z$, and $\sigma$ is globally Lipschitz in $z$; 
that $b, b_\mathrm{I}, \sigma$ satisfy the growth conditions
\begin{align*}		
\left| b(x, z) \right| + \left| b_\mathrm{I}(x, z) \right| + \left| \sigma \sigma^* (x, z) \right| &\leq C (1 + |z|)^\beta, \\
\sum_{|k| = 1}^2 \left| D^k_x b(x, z) \right| +\left| D^k_x \sigma \sigma^* (x, z) \right| &\leq C (1 + |z|^q), \\
\sum_{|k| = 1}^3 \left| D^k_x b_\mathrm{I}(x, z) \right| &\leq C (1 + |z|^q),
\end{align*}
for some \Nicolas{$\beta < -2$} and $q > 0$;
that $h$ is bounded in $(x, z)$, $h(\cdot, z) \in C^3$ for each $z$, and $h$ is globally Lipschitz in $z$.
\item $\overbar{a} + \wt{a} \succ 0$ uniformly in $x$;
\hyperlink{H}{$H^{2, 2 + \alpha}$} holds for $\alpha \in (0, 1)$;
for each $z$, $b(\cdot, z), b_\mathrm{I}(\cdot, z), \sigma(\cdot, z) \in C^2$; 
that $b$ and $b_\mathrm{I}$ are Lipschitz in $z$, and $\sigma$ is globally Lipschitz in $z$;
that $b, b_\mathrm{I}, \sigma$ satisfy the growth conditions
\begin{align*}		
\left| b(x, z) \right| + \left| b_\mathrm{I}(x, z) \right| + \left| \sigma \sigma^* (x, z) \right| &\leq C (1 + |z|)^\beta, \\
\sum_{|k| = 1}^2 \left| D^k_x b(x, z) \right| + \left| D^k_x b_\mathrm{I}(x, z) \right| + \left| D^k_x \sigma \sigma^* (x, z) \right| &\leq C (1 + |z|^q),
\end{align*}
for some \Nicolas{$\beta < -2$} and $q > 0$;
$h$ is bounded in $(x, z)$, that $h$ is globally Lipschitz in $(x, z)$. 
If $a \succ 0$, which implies $\overbar{a} + \wt{a} \succ 0$, then the Lipschitz condition in $z$ for $b, b_\mathrm{I}$ can be relaxed to $\alpha$-H\"older continuity. 
\end{enumerate}
Then there exists a metric $d$ on $C([0, T]; P(\R{m}))$, the space of continuous processes from $[0, T]$ to the space of probability measures on $\R{m}$, that generates the topology of weak convergence, such that $\pi^{\epsilon, x} \rightarrow \pi^0$ in probability. 
\begin{proof}
From Theorem \ref{theorem: weak convergence of unnormalized filter} and Lemma \ref{lemma: weak convergence of unnormalized filter implies convergence of normalized} we get $\pi^{\epsilon, x} - \pi^0 \Rightarrow 0$ as $\epsilon \rightarrow 0$.
$\pi^{\epsilon, x}$ and $\pi^0$ are random variables in the space $C([0, T]; P(\R{m}))$. 
We define a continuous bounded metric $d$ on this space that generates the topology of weak convergence as follows, 
\begin{align*}
d(\mu, \nu) = 1 \wedge ( \sup_{0 \leq t \leq T} \wt{d}(\mu_t, \nu_t) ),
\end{align*}
where we assume that $\wt{d}$ is a translation invariant metric that generates the topology of weak convergence on $P(\R{m})$.
Then $d$ inherits the translation invariant property from $\wt{d}$, and from weak convergence of $\pi^{\epsilon, x} - \pi^0$ in the space of signed measures, we have
\[
\lim_{\epsilon \rightarrow 0 } \E[\Q{}] { d( \pi^{\epsilon, x} , \pi^0 ) } = \lim_{\epsilon \rightarrow 0 } \E[\Q{}] { d( \pi^{\epsilon, x} - \pi^0 , 0 ) } = 0.
\]
And therefore we retrieve convergence in probability,
\begin{align*}
\lim_{\epsilon \rightarrow 0} \Q{} \left( d( \pi^{\epsilon, x} , \pi^0 ) \geq \delta \right) 
\leq \frac{1}{\delta} \lim_{\epsilon \rightarrow 0 } \E[\Q{}] { d( \pi^{\epsilon, x} , \pi^0 ) } 
= 0, \qquad \text{for each $\delta > 0$}. 
\end{align*}
It remains to show that there is a translation invariant metric $\wt{d}$ that generates the weak topology on $P(\R{m})$. 
For this, we can borrow the argument from \cite[Corollary 6.9, p.2322]{Imkeller:2013}. 
\end{proof}
}
}

\subsection{Change of Probability Measure and Zakai Equation}

To prove the main result, the analysis which we perform will actually be concerned with unnormalized conditional measures that are defined using change of probability measure transformations. 
The new collection of measures on the filtered probability space are denoted by $(\Pe)$. For any fixed $\epsilon$, $\Pe$ and $\Q{}$ will be mutually absolutely continuous with Radon-Nikodym derivatives 
\begin{align*}
D^\eps_t 
\equiv \restrict{\frac{d\P^\eps}{d\Q{}}}{\mF_t} 
= \exp \left( - \int_0^t \langle h(X^\eps_s, Z^\eps_s), dB_s \rangle - \frac{1}{2} \int_0^t \left| h(X^\eps_s, Z^\eps_s) \right|^2 ds \right),
\end{align*}

\begin{align*}
\wt{D}^\eps_t 
\equiv (D^\eps_t)^{-1} 
= \restrict{\frac{d\Q{}}{d\P^\eps}}{\mF_t} 
= \exp \left( \int_0^t \langle h (\Xe_s, \Ze_s), d\Ye_s \rangle - \frac{1}{2} \int_0^t \left| h (X^\eps_s, Z^\eps_s) \right|^2 ds \right). \\
\end{align*}
Then by Girsanov's theorem, under $\P^\eps$ the process $\Ye$ is a Brownian motion.
For a fixed test function $\varphi \in C^2_b(\R{m} \times \R{n}; \R{})$ and time $t \in [0 , T]$, we characterize the unnormalized conditional measure $\rho^\epsilon_t$ as, 
\begin{align*}
\rho^\eps_t(\varphi) 
= \E[\Pe]{\varphi(\Xe_t, \Ze_t)\wt{D}^\eps_t \st \mYe_t},
\end{align*}
and its relation to $\pi^\epsilon$ through the Kallianpur-Striebel formula, 
\begin{align*}\label{equation: Kallianpur-Striebel}
\pi^\eps_t(\varphi) 
= \frac{\E[\Pe]{\varphi(\Xe_t, \Ze_t)\wt{D}^\eps_t \st \mYe_t}}{\E[\Pe]{ \wt{D}^\eps_t \st \mYe_t}} 
= \frac{\rho^\eps_t(\varphi)}{\rho^\eps_t(1)},
\qquad \forall t \in [0, T], 
\qquad \text{$\Q{}, \Pe$-a.s}.
\end{align*}
The action of $\rho^\epsilon$ on test functions $\varphi \in C^2_b(\R{m} \by \R{n}; \R{})$ gives the Zakai evolution equation, 
\EquationAligned{
\label{equation: Zakai}
\rho^\epsilon_t(\varphi) &= \rho^\epsilon_0(\varphi) + \int_0^t \rho^\epsilon_s \left( \mG^\epsilon \varphi \right) ds + \int_0^t \langle \rho^\epsilon_s ( \varphi h + \alpha \sigma^* \nabla_x \varphi ),  d\Ye_s \rangle, \\
\rho^\epsilon_0(\varphi) &= \E[\Q{}]{ \varphi(\Xe_0,\Ze_0) }.
}
When $\varphi \in C^2_b(\R{m}; \R{})$, we consider the $x$-marginal, 
\begin{align*}
\rho^{\epsilon, x}_t(\varphi) = \int \varphi(x) \rho^\epsilon_t(dx, dz),
\end{align*}
which is related to $\pi^{\epsilon, x}$ through the Kallianpur-Striebel formula, 
\begin{align*}
\pi^{\epsilon, x}_t(\varphi) 
= \frac{\rho^{\epsilon, x}_t(\varphi)}{\rho^{\epsilon, x}_t(1)},
\qquad \forall t \in [0, T], 
\qquad \text{$\Q{}, \Pe$-a.s}.
\end{align*}
This is easy to see since $\rho^{\epsilon}(1)  = \rho^{\epsilon, x}(1)$.

We next define the averaged unnormalized filter $\rho^0$ as the solution of the following evolution equation, 
\EquationAligned{
\label{equation: homogenized Zakai}
\rho^0_t(\varphi) &= \rho^0_0(\varphi) + \int_0^t \rho^0_s(\mG^\dagger \varphi) ds + \int_0^t \langle \rho^0_s(\varphi \overbar{h} + \alpha \overbar{\sigma}^* \nabla_x \varphi ),  d\Ye_s \rangle, \\
\rho^0_0(\varphi) &= \E[\Q{}]{ \varphi(X^0_0) },
}
where $\varphi \in C^2_b(\R{m}; \R{})$.
And then by the Kallianpur-Striebel formula we relate the averaged (normalized) filter $\pi^0$ to the unnormalized variant, 
\begin{align*}
\pi^0_t(\varphi) = \frac{\rho^0_t(\varphi) }{\rho^0_t(1) },
\qquad \forall t \in [0, T], 
\qquad \text{$\Q{}, \Pe$-a.s}.
\end{align*}
The uniqueness of $\rho^0$ follows from the same assumptions and proof to be given in Lemma \ref{lemma: uniqueness of weak limits}. 

We will later show in Lemma \ref{lemma: weak convergence of unnormalized filter implies convergence of normalized} that under appropriate assumptions, weak convergence of $\rho^{\epsilon, x} - \rho^0$ to zero will imply weak convergence of $\pi^{\epsilon, x} - \pi^0$ to zero, and therefore we can focus on showing convergence of the unnormalized difference for the main analysis.

\subsubsection{Representation of the Averaged Unnormalized Conditional Distribution}
\label{section: representation of averaged unnormalized conditional distribution (chapter: quantitative)}

Just as $\pi^0$ is not the filter for the averaged system, $\rho^0$ is not the unnormalized conditional distribution for the averaged system, and therefore a representation of this measure acting on $\varphi \in C^2_b$ test functions as conditional expectation requires a bit more work.
But such a representation will be necessary for the computation of some estimates. 

To get such a representation, we introduce a signal process $X^0$ to be a diffusion process with infinitesimal generator $\mG^\dagger$. 
Therefore, consider the following SDE, 
\begin{align}
\label{equation: SDE to achieve averaged unnormalized filter equation}
dX^0_t 
&= \left[ \overbar{b}(X^0_t) + \wt{b}(X^0_t) \right] dt + \wt{a}^{1/2}(X^0_t) d\widetilde{W}_t + ( \overbar{a}(X^0_t) - \overbar{\sigma} \overbar{\sigma}^*(X^0_t) )^{1/2} d\widehat{W}_t + \overbar{\sigma}(X^0_t) dW_t, \\
X^0_0 &\sim \Q{}_{\Xe_0}. \nonumber
\end{align}
Here $\widetilde{W}$ and $\widehat{W}$ are new $m$-dimensional independent Brownian motions, independent of  $(V, W, U)$ under $\Q{}$ as well as independent of the initial condition $\Q{}_{\Xe_0}$.
The Cholesky factor $( \overbar{a}(X^0_t) - \overbar{\sigma} \overbar{\sigma}^*(X^0_t) )^{1/2}$ exists, since from an application of Jensen's inequality, one can show that $\overbar{a}(x) - \overbar{\sigma} \overbar{\sigma}^*(x)$ is positive semidefinite for each $x \in \R{m}$. 
For this dispersion coefficient to be Lipschitz continuous, we require $(\overbar{a} - \overbar{\sigma} \overbar{\sigma}^*) \in C^2_b$ \cite[Lemma 2.3.3]{Stroock:2008c}. 
This will be true if $\sigma(\cdot, z) \in C^2_b$ for each $z$, which is assumed in Theorem \ref{theorem: main result, weak convergence of the normalized filter}.

{
\rem{
An interesting observation regarding Eq. \ref{equation: SDE to achieve averaged unnormalized filter equation}, is that we may have $\overbar{\sigma} = 0$, and this implies that the SDE for the averaged filter may have no correlation at all, or less correlation than the original system. 
}
}

We now define the process
\begin{align*}
\wt{D}^0_t 
= \exp \left( \int_0^t \langle \overbar{h} (X^0_s), d\Ye_s \rangle - \frac{1}{2} \int_0^t \left| \overbar{h} (X^0_s) \right|^2 ds \right),
\end{align*}
which is used to give the representation of $\rho^0$ on $C^2_b$ test functions as follows, 
\begin{align*}
\rho^0_t(\varphi) 
= \E[\Pe]{\varphi(X^0_t)\wt{D}^0_t \st \mYe_t}.
\end{align*}

\section{Preliminary Estimates}
\label{section: preliminary estimates}

In this section we provide several preliminary estimates that will be needed for the main analysis. 
Some additional comments regarding notation are first introduced and then some assumptions are defined. 

The relation $a \lesssim b$ will indicate that $a \leq C b$ for a constant $C > 0$ that is independent of $a$ and $b$, but that may depend on parameters that are not critical for the bound being computed.
We will use the notation $T^{F, x}$ for the semigroup of $Z^x$, and denote processes with $Z^{\epsilon, x; (s, z)}$ to represent the process $Z^{\epsilon, x}$ started at time $s$ at $z \in \R{n}$. 
We will say that a function $\theta(x, z)$ is centered with respect to $\mu_\infty$ (the family of invariant measures parameterized by $x \in \R{m}$) if 
\[
\int_{\R{n}} \theta(x, z) \mu(dz; x) = 0, \qquad \forall x \in \R{m}. 
\]

Let $H_L$ denote the assumption that for each $K > 0$, there exists a constant $C_K$ such that for all $x, x' \in \R{m}$, $|z| \leq K$: 
\begin{align}
\tag{$H_L$}
\label{assumption: local Lipschitz of slow coefficients}
|b(x, z) - b(x', z)| + |b_\mathrm{I}(x, z) - b_\mathrm{I}(x', z)| + |\sigma(x, z) - \sigma(x', z)| \leq C_K |x - x'|. 
\end{align}

Let $H_P$ denote the assumption that there exists $K, \alpha, p_1, p_2 > 0$ such that for all $(x, z) \in \R{m} \times \R{n}$: 
\begin{align}
\tag{$H_P$}
\label{assumption: linear and polynomial growth of slow drift and dispersion coefficients}
|b(x, z)| &\leq K (1 + |x|) (1 + |z|^{p_1}), \\
|\sigma(x, z)| &= \sqrt{ \op{Tr}(\sigma \sigma^*(x, z)) } \leq K (1 + |x|^{1/2}) (1 + |z|^{p_2}) . \nonumber
\end{align} \ifbool{ShowFootnotes}{ \footnote{Because we don't have correlation between slow and fast, we may be able to get $(1 + |x|)$} }{}
Note that from $H_P$ we have $|\sigma(x, z)| \lesssim 1 + |x| + |z|^{2 p_2}$ and hence implies a linear growth in $x$ and polynomial growth in $z$. 
Also from $H_P$, $| \sigma \sigma^*(x, z)| \lesssim (1 + |x|^2 + |z|^{4p_2})$. 

Let $H_I$ denote the assumption that for some $K, p > 0$, $b_\mathrm{I}$ satisfies the following growth condition, 
\begin{align}
\tag{$H_I$}
\label{assumption: growth condition for intermediate coefficient}
\sum_{|\alpha| \leq 2} \sup_{x \in \R{m}} | D^\alpha_x b_\mathrm{I}(x, z) | \leq K (1 + |z|^p ).
\end{align}

The next result is from \cite[p.1172]{Pardoux:2003} and provides the result that $\Xe \Rightarrow X^0$ in the limit $\epsilon \rightarrow 0$. 
{
\thm{
\label{theorem: pardoux and veretennikov diffusion approximation}
Let $(\Xe, \Ze)$ satisfy the stochastic differential equations of Eq. \ref{equation: multiscale signal process} with initial conditions $(\Xe_0, \Ze_0) = (x, z) \in \R{m} \times \R{n}$ for each $\epsilon \in (0, 1)$.
Assume \ref{assumption: positive recurrence}, \ref{assumption: uniform ellipticity}, \hyperlink{H}{$H^{2, 2 + \alpha}$} for $\alpha \in (0, 1)$, \ref{assumption: local Lipschitz of slow coefficients}, and \ref{assumption: linear and polynomial growth of slow drift and dispersion coefficients}. 
Let $b_\mathrm{I} \in C^{2, \alpha}$ satisfy \ref{assumption: growth condition for intermediate coefficient} and be centered with respect to $\mu_\infty$. 
Then for any $T > 0$, the process $\Xe$ converges weakly in the limit $\epsilon \rightarrow 0$, to the Markov process $X^0$ with generator $\mG^\dagger$. 

\begin{proof}
See remarks in Section \ref{section: remark on condition for the fast semigroup} and \cite[p.1172]{Pardoux:2003}. 
\end{proof}
}
}

\subsection{Estimates with the Fast Semigroup}

{
\lemma{
\label{lemma: bounds on functions under semigroup}
Assume \hyperlink{H}{$H^{k, l}$}, with $k \in \N{}_0, l \in \N{}$, and let $\theta \in C^{k, j}(\R{m} \times \R{n}; \R{})$ for $j \leq l$ satisfy for some $C, p > 0$
\[
\sum_{|\alpha| \leq k} \sum_{|\beta| \leq j} | D^\alpha_x D^\beta_z \theta(x, z) | \leq C (1 + |x|^p + |z|^p).
\]
Then 
\[
(t, x, z) \mapsto T^{F, x}_t \left( \theta (x, \cdot) \right) (z) \in C^{0, k, j}(\R{+} \times \R{m} \times \R{n}; \R{})
\]
and there exist $C_1, p_1 > 0$, such that for all $(t, x, z) \in [0, \infty) \times \R{m} \times \R{n}$
\[
\sum_{|\alpha| \leq k} \sum_{|\beta| \leq j} | D^\alpha_x D^\beta_z T^{F, x}_t \left( \theta (x, \cdot) \right) (z) | \leq C_1 e^{C_1 t} (1 + |x|^{p_1} + |z|^{p_1}).
\]
If the bound on the derivatives of $\theta$ can be chosen uniformly in $x$, that is, 
\[
\sum_{|\alpha| \leq k} \sum_{|\beta| \leq j} \sup_x | D^\alpha_x D^\beta_z \theta(x, z) | \leq C (1 + |z|^p),
\]
then the bound on the derivatives of $T^{F, x}_t \left( \theta (x, \cdot) \right) (z)$ is also uniform in $x$, 
\[
\sum_{|\alpha| \leq k} \sum_{|\beta| \leq j} \sup_x | D^\alpha_x D^\beta_z T^{F, x}_t \left( \theta (x, \cdot) \right) (z) | \leq C_1 e^{C_1 t} (1 + |z|^{p_1}).
\]

\begin{proof}
The proposition is a slight generalization of \cite[Proposition 5.1]{Imkeller:2013}. 
The proof is the same as in \cite[Proposition 5.1]{Imkeller:2013}.
\end{proof}
}
}

{
\lemma{
\label{lemma: regularity and bound on averaged function}
Assume \ref{assumption: positive recurrence}, \ref{assumption: uniform ellipticity} and \hyperlink{H}{$H^{k, 2 + \alpha}$} for $\alpha \in (0, 1)$, and $k \in \N{}_0$. 
Let $\theta \in C^{k, 0}(\R{m} \times \R{n}; \R{})$ satisfy for some $C, p > 0$,
\[
\sum_{|\gamma| \leq k} \sup_x | D^\gamma_x \theta(x, z) | \leq C ( 1 + |z|^p ). 
\]
Then 
\[
x \mapsto \mu_\infty(\theta; x)(x') = \int_{\R{n}} \theta(x', z) \mu_\infty(dz; x) = \int_{\R{n}} \theta(x', z) p_\infty(z; x) dz \in C^{k}_b(\R{m}; \R{}).
\]

\begin{proof}
The proposition is a slight generalization of a part of \cite[Proposition 5.2]{Imkeller:2013} and the proof follows the same argument as given there.
\end{proof}
}
}

{
\lemma{
\label{lemma: bounds on centered functions under semigroup}
Assume \ref{assumption: positive recurrence}, \ref{assumption: uniform ellipticity} and \hyperlink{H}{$H^{k, 2 + \alpha}$} for $\alpha \in (0, 1)$ and $k \in \N{}_0$. 
Let $j \in \{ 0, 1 \}$, and $\theta \in C^{k, j + \alpha (1 - j)}(\R{m} \times \R{n}; \R{})$ satisfy the growth condition, 
\[
\sum_{|\gamma| \leq k} \sum_{|\beta| \leq j} \sup_x | D^\gamma_x D^\beta_z \theta(x, z) | \leq C (1 + |z|^p),
\]
for some $C, p > 0$ .
Assume additionally that $\theta$ satisfies the centering condition, 
\[
\int_{\R{n}} \theta(x, z) \mu_\infty(dz; x) = 0, \quad \forall x \in \R{m}.
\]
Then 
\[
(x, z) \mapsto \int_0^\infty T^{F, x}_t(\theta(x, \cdot))(z) dt \in C^{k, j}(\R{m} \times \R{n}; \R{}),
\]
and for every $q > 0$ there exist $C', q' > 0$, such that, 
\[
\sum_{|\gamma| \leq k} \sum_{|\beta| \leq j} \int_0^\infty \sup_x | D^\gamma_x D^\beta_z T^{F, x}_t(\theta(x, \cdot))(z) |^q dt \leq C' (1 + |z|^{q'}).
\]
\begin{proof}
The proposition is a slight generalization of a part of \cite[Proposition 5.2]{Imkeller:2013} and the proof follows the same argument as given there.
\end{proof}
}
}

\subsection{Estimates on SDE Solutions}

{
\lemma{
\label{lemma: moments on fast process}
Assume $f$ is bounded and that $f$ and $gg^*$ are H\"older continuous in $z$ uniformly in $x$ for some uniform constant.
Assume that the conditions \ref{assumption: positive recurrence} and \ref{assumption: uniform ellipticity} hold. 
Then for any $p > 0$ there exists $C_p > 0$ such that 
\begin{align*}
\sup_{(t, \epsilon, x) \in [0, \infty) \times (0, 1] \times \R{m}} \E{ | \Ze_t |^p}[(\Xe_0, \Ze_0) = (x, z)] \lesssim 1 + |z|^p.
\end{align*}
\begin{proof}
The proposition is a slight generalization of a part of \cite[Proposition 5.3]{Imkeller:2013} and the proof follows the same argument as given there.
\end{proof}
}
}

{
\lemma{
\label{lemma: moment on intermediate scale forcing}
Assume the conditions \ref{assumption: positive recurrence}, \ref{assumption: uniform ellipticity} and \hyperlink{H}{$H^{2, 2 + \alpha}$} for some $\alpha \in (0, 1)$; that $b, \sigma$ are bounded for all $(x, z)$; and that $b_\mathrm{I} \in C^{2, 1}(\R{m} \times \R{n}; \R{m})$ that satisfies the centering condition, 
\begin{align*}
\int_{\R{n}} b_\mathrm{I}(x, z) \mu_\infty(dz; x) = 0,
\end{align*}
where $\mu_\infty(x)$ is the unique stationary distribution for the process $Z^x$, and that for some $C, q_1 > 0$, it has the following growth condition, 
\begin{align*}
\sum_{|\alpha| \leq 2} \sum_{|\beta| \leq 1} \sup_x | D^\alpha_x D^\beta_z b_\mathrm{I}(x, z) | \leq C (1 + |z|^{q_1}).
\end{align*}
Then for every $p \geq 2$ there exists $q > 0$ such that for $0 \leq r < t < \infty$ 
\begin{align*}
\E* \left| \frac{1}{\epsilon} \int_r^t b_\mathrm{I}(\Xe_s, \Ze_s) ds  \right|^p 
\lesssim \epsilon^p (1 + |z|^q) 
&+ (t - r)^{p - 1}  (1 + \epsilon^p) \int_r^t 1 + \E* | \Ze_s |^q ds \\
&+ (t - r)^{(p/2) - 1} (1 + \epsilon^p) \int_r^t 1 + \E* | \Ze_s |^q ds. 
\end{align*}
\begin{proof}
We start by considering the solution of the following backward partial differential equation, 
\begin{align*}
- \partial_s \psi_s(x, z) = \frac{1}{\epsilon^2} \mG_F \psi_s(x, z) + \frac{1}{\epsilon} b_\mathrm{I}(x, z), \quad \psi_t(x, z) = 0.
\end{align*}
The solution of which is given by a Feynman-Kac representation, 
\FutureWork{why not use Poisson equation here, seems easier?}
\begin{align*}
\psi_s(x, z) 
=  \E* \int_s^t \frac{1}{\epsilon} b_\mathrm{I}(x, Z^{\epsilon, x; (s, z)}_r) dr 
&= \frac{1}{\epsilon} \int_s^t T_{(r - s) / \epsilon^2}^{F, x} \left( b_\mathrm{I}(x, \cdot) \right)(z) dr 
= \epsilon \int_0^{(t - s) / \epsilon^2} T^{F, x}_u \left( b_\mathrm{I}(x, \cdot) \right)(z) du,
\end{align*}
where $T^{F, x}$ is the semigroup associated with the process $Z^x$ and a change of time has been used in the last equality relation. 
Then since $b_\mathrm{I}$ satisfies the conditions of Lemma \ref{lemma: bounds on centered functions under semigroup}, we have 
\begin{align*}
\sum_{|\alpha| \leq 2} \sum_{|\beta| \leq 1} \sup_{s \in [0, t]} | D^\alpha_x D^\beta_z \psi_s(x, z) |^p \lesssim \epsilon^p (1 + |z|^{q_2} ),
\end{align*}
for some $q_2 > 0$. 
Applying It\^o's formula to $\psi_t(x, z)$ gives, 
\begin{align*}
0 = \psi_r(x, z) &+ \frac{1}{\epsilon^2} \int_r^t \mG_F \psi_s(\Xe_s, \Ze_s) ds + \frac{1}{\epsilon} \int_r^t \nabla_z \psi_s(\Xe_s, \Ze_s)  g(\Xe_s, \Ze_s) dV_s \\
&+ \int_r^t \mG^\epsilon_S \psi_s(\Xe_s, \Ze_s) ds + \int_r^t \nabla_x \psi_s(\Xe_s, \Ze_s) \sigma(\Xe_s, \Ze_s) dW_s  \\
&- \frac{1}{\epsilon^2} \int_r^t \mG_F \psi_s(\Xe_s, \Ze_s) ds - \frac{1}{\epsilon} \int_r^t b_\mathrm{I}(\Xe_s, \Ze_s) ds.
\end{align*}
Eliminating terms and rearranging simplifies to 
\begin{align*}
\frac{1}{\epsilon} \int_r^t b_\mathrm{I}(\Xe_s, \Ze_s) ds 
= \psi_r(x, z) &+ \int_r^t \mG_S \psi_s(\Xe_s, \Ze_s) ds + \frac{1}{\epsilon} \int_r^t \nabla_x \psi_s(\Xe_s, \Ze_s) b_\mathrm{I}(\Xe_s, \Ze_s) ds \\
&+ \frac{1}{\epsilon} \int_r^t \nabla_z \psi_s(\Xe_s, \Ze_s)  g(\Xe_s, \Ze_s) dV_s  + \int_r^t \nabla_x \psi_s(\Xe_s, \Ze_s) \sigma(\Xe_s, \Ze_s) dW_s .
\end{align*}
The first term will contribute, 
\begin{align*}
\E* | \psi_r(x, z) |^p \leq \sup_{s \in [0, t]} | \psi_s(x, z) |^p \lesssim \epsilon^p (1 + |z|^{q_2}).
\end{align*}
From the boundedness of $b$ and $\sigma$ we have for the second term
\begin{align*}
\E* \left| \int_r^t \mG_S \psi_s(\Xe_s, \Ze_s) ds \right|^p 
&\leq (t - r)^{p - 1} | b |_\infty^p \int_r^t \E* \left| \nabla_x \psi_s(\Xe_s, \Ze_s) \right|^p ds + (t - r)^{p - 1} | \sigma |_\infty^p \int_r^t \E* \left| \nabla^2_x \psi_s(\Xe_s, \Ze_s) \right|^p ds \\
&\lesssim (t - r)^{p - 1} ( | b |_\infty^p + | \sigma |_\infty^p ) \epsilon^p \int_r^t 1 + \E* \left| \Ze_s \right|^{q_2} ds.
\end{align*}
For the third term, 
\begin{align*}
\left| \nabla_x \psi_s (\Xe_s, \Ze_s) b_\mathrm{I}(\Xe_s, \Ze_s) \right|^p 
\leq | \nabla_x \psi_s(\Xe_s, \Ze_s) |^p | b_\mathrm{I}(\Xe_s, \Ze_s)  |^p 
\lesssim \epsilon^p \left( 1 + | \Ze_s |^{q_2} \right) \left( 1 + | \Ze_s |^{q_1} \right)^p
\end{align*}
and therefore 
\begin{align*}
\E* \left| \frac{1}{\epsilon} \int_r^t \nabla_x \psi_s(\Xe_s, \Ze_s) b_\mathrm{I}(\Xe_s, \Ze_s) ds \right|^p 
\lesssim (t - r)^{p - 1} \int_r^t 1 + \E* \left| \Ze_s \right|^{q_3} ds.
\end{align*}
The stochastic integrals follow in a similar manner after application of the Burkholder-Davis-Gundy (BDG) inequality, and the boundedness of $\sigma$ and $g$, 
\begin{align*}
\E* \left| \int_r^t \nabla_x \psi_s(\Xe_s, \Ze_s) \sigma(\Xe_s, \Ze_s) dW_s \right|^p
&= \E* \left( \sup_{u \in [r, t]} \left| \int_r^u \nabla_x \psi_s(\Xe_s, \Ze_s) \sigma(\Xe_s, \Ze_s) dW_s \right| \right)^p \\
&\leq C_p \E* \left( \int_r^t \left| \nabla_x \psi_s(\Xe_s, \Ze_s) \sigma(\Xe_s, \Ze_s) \right|^2 ds \right)^{p / 2} \\
&\leq C_p (t - r)^{(p/2) - 1} \epsilon^p | \sigma |_\infty^p  \int_r^t 1 + \E* | \Ze_s |^{q_2} ds .
\end{align*}
The bound for the other stochastic integral follows in the same manner, 
\begin{align*}
\E* \left| \frac{1}{\epsilon} \int_r^t \nabla_z \psi_s(\Xe_s, \Ze_s) g(\Xe_s, \Ze_s) dV_s \right|^p 
\leq C_p (t - r)^{(p / 2) - 1} | g |_\infty^p \int_r^t 1 + \E* | \Ze_s |^{q_2} ds. 
\end{align*}
Collecting all the terms, now yields the desired result.
\end{proof}
}
}

{
\lemma{
\label{lemma: moments on slow process}
Assume the same setup as Lemma \ref{lemma: moment on intermediate scale forcing},
then for every $p \geq 2$ there exists $q > 0$ such that for $T > 0$ 
\begin{align*}
\sup_{(t, \epsilon) \in [0, T] \times (0, 1]} \E{ | \Xe_t |^p}[(\Xe_0, \Ze_0) = (x, z)] \lesssim 1 + |x|^p + |z|^q. 
\end{align*}
\begin{proof}
Since $b$ and $\sigma$ are bounded, we get, 
\begin{align*}
\E* | \Xe_t |^p \lesssim 1 + \E* | \Xe_0 |^p + \E* \left| \int_0^t \frac{1}{\epsilon} b_\mathrm{I}(\Xe_s, \Ze_s) ds \right|^p.
\end{align*}
Using the result of Lemma \ref{lemma: moment on intermediate scale forcing} for the moment of the intermediate scale forcing and then Lemma \ref{lemma: moments on fast process} for the moment of the fast process, we get
\begin{align*}
\E* | \Xe_t |^p 
&\lesssim \E* | \Xe_0 |^p + \epsilon^p (1 + \E* | \Ze_0 |^q) + (1 + \epsilon^p) \int_0^t 1 + \E* | \Ze_s |^q ds 
\lesssim \E* | \Xe_0 |^p + (1 + \epsilon^p) (1 + \E* | \Ze_0 |^q).
\end{align*}
Repeating the proof, but conditioning on $(\Xe_0, \Ze_0) = (x, z)$ gives the desired result. 
\end{proof}
}
}

{
\lemma{
\label{lemma: moments on increments of slow process}
Assume the same setup as Lemma \ref{lemma: moment on intermediate scale forcing}, then for $|t - s| \leq 1$ and $p \geq 2$, there exists a $q \geq 0$ such that
\begin{align*}
\E{ | \Xe_t - \Xe_s |^p}[(\Xe_0, \Ze_0) = (x, z)] \lesssim \epsilon^p ( 1 + |z|^q ) + (t - s)^{p/2} (1 + \epsilon^p) (1 + |z|^q).
\end{align*}
\begin{proof}
Without loss of generality, assume $s < t$.
\begin{align*}
\E{ | \Xe_t - \Xe_s |^p }
&\lesssim (t - s)^{p - 1} \int_s^t \E*{} \left| b(\Xe_u, \Ze_u) \right|^p du + \E*{} \left| \frac{1}{\epsilon} \int_s^t b_\mathrm{I}(\Xe_u, \Ze_u) du \right|^p +  \E*{} \left| \int_s^t \sigma(\Xe_u, \Ze_u) dW_u \right|^p \\
&\leq (t - s)^p | b |_\infty^p + (t - s)^{p / 2} | \sigma |_\infty^p + \E*{} \left| \frac{1}{\epsilon} \int_s^t b_\mathrm{I}(\Xe_u, \Ze_u) du \right|^p. 
\end{align*}
Now using Lemma \ref{lemma: moment on intermediate scale forcing}, Lemma \ref{lemma: moments on fast process}, we have 
\begin{multline*}
\E{ | \Xe_t - \Xe_s |^p }
\lesssim (t - s)^p + (t - s)^{p / 2} + \epsilon^p ( 1 + \E* | \Ze_0 |^q ) + (t - s)^p (1 + \epsilon^p) (1 + \E* | \Ze_0 |^q ) + (t - s)^{p/2} (1 + \epsilon^p) (1 + \E* | \Ze_0 |^q ) \\
\lesssim \epsilon^p ( 1 + \E* | \Ze_0 |^q ) + (t - s)^{p/2} (1 + \epsilon^p) (1 + \E* | \Ze_0 |^q ).
\end{multline*}
Repeating the proof, but conditioning on $(\Xe_0, \Ze_0) = (x, z)$ gives the desired result. 
\end{proof}
}
}

{
\lemma{
\label{lemma: moments of Girsanov} 
Assume $h$ is bounded, then for $p \geq 2$ and $T > 0$, 
\begin{align*}
\sup_{\epsilon \in (0, 1]} \sup_{t \leq T} \E*[\Pe] \left| \wt{D}^\epsilon_t \right|^p < \infty
\quad \mathrm{and} \quad 
\sup_{t \leq T} \E*[\Pe]{ \left| \wt{D}^0_t \right|^p } < \infty.
\end{align*}
Further, for $| t - s | < 1$, we have 
\begin{align*}
\sup_{\epsilon \in (0, 1]}
\E*[\Pe] \left| \wt{D}^\epsilon_t - \wDe_s \right|^p 
\lesssim C_p (t - s)^{p / 2} | h |_\infty^p < \infty.
\end{align*}
\begin{proof}
We have that $\wt{D}^\epsilon_t$ satisfies 
\begin{align*}
\wDe_t = 1 + \int_0^t \wDe_s \langle h(\Xe_s, \Ze_s), d\Ye_s \rangle.
\end{align*}
Using the boundedness of $h$, the first result now follows from an application of BDG and Gr\"onwall's lemma. 
The same proof applies for the moment bound of $\wt{D}^0$.
The bound for the increment now follows, 
\begin{align*}
\E*[\Pe] \left| \wt{D}^\epsilon_t - \wDe_s \right|^p 
&= \E*[\Pe] \left| \int_s^t \wDe_u \langle h(\Xe_u, \Ze_u), d\Ye_u \rangle \right|^p \\
&\leq C_p (t - s)^{(p/2) - 1} \int_s^t \E[\Pe]{ \left| \wDe_u \right|^p \left| h(\Xe_u, \Ze_u) \right|^p } du 
\lesssim  C_p (t - s)^{p / 2} | h |_\infty^p < \infty.
\end{align*} 
\FutureWork{
\Nicolas{
The case $1 \leq p < 2$ follows from $L^p(\Pe) \subset L^2(\Pe)$. 
}
}
\end{proof}
}
}

\subsection{Estimates using the Poisson Equation}

{
\thm{
\label{theorem: solution, bounds and regularity of Poisson equation}
Consider the Poisson equation, 
\begin{align*}
\mG u(x, z) = - \psi(x, z), 
\end{align*}
where $x \in \R{m}$ is a parameter, $\mG$ is the generator
\begin{align*}
\mG(x, z) \equiv \sum_{i=1}^n f_i(x, z) \frac{\partial}{\partial z_i} + \frac{1}{2} \sum_{i, j=1}^n (gg^*)_{ij}(x, z) \frac{\partial^2}{\partial z_i \partial z_j},
\end{align*}
and $\psi \in C^{k, \alpha}$ for $k \geq 1$ and $\alpha > 0$. 
Assume that $f, g$ satisfy the assumptions of \ref{assumption: positive recurrence}, \ref{assumption: uniform ellipticity} and \hyperlink{H}{$H^{1 , 2 + \alpha}$}. 
Let $\mu_\infty(x)$ be the unique stationary distribution of $Z^x$ for each fixed $x \in \R{m}$,
and assume that $\psi$ is centered for each $x \in \R{m}$, 
\begin{align*}
\int_{\R{n}} \psi(x, z) \mu_\infty(dz; x) = 0.
\end{align*}
Further, assume the growth conditions
\begin{align*}		
\left| \psi(x, z) \right| &\leq C_0 (1 + |z|)^{\Nicolas{\beta}}, \\
\sum_{|j| = 1}^k \left| D^j_x \psi(x, z) \right| &\leq C_1(x) (1 + |z|^q),
\end{align*}
for some \Nicolas{$\beta < -2$}\ifbool{ShowFootnotes}{\footnote{$\beta < -2$ is given in \cite{Pardoux:2003}, though it seems like it may be possible with something weaker e.g. $\beta < 0$ based on \cite{Pardoux:2001vc}}}{} and $q > 0$. 
Then the solution of the Poisson equation exists, belongs to the Sobolev space $\cap_{p \in (1, \infty)} W^2_{p, \mathrm{loc}}$, is unique up to an additive constant such that for any $x$ the centering condition 
\begin{align*}
\int_{\R{n}} u(x, z) \mu_\infty(dz; x) = 0
\end{align*}
holds, the solution satisfies $u(\cdot, z) \in C^k$ for any $z$, and the following holds true for some $q', q''$ and some constants $C_2, C_3(x)$, 
\begin{align*}
\Nicolas{|u(x, z)|} &\leq C_2, \\
\sum_{|j| = 1}^k |D^j_x u(x, z)| &\leq C_3(x) (1 + |z|^{q'}), \\
|\nabla_x \nabla_z u(x, z)| &\leq C_3(x) (1 + |z|^{q''}).
\end{align*}
\begin{proof} 
This is a combination of Proposition 1 and Theorem 3 from \cite{Pardoux:2003}, but restricted for our needs. 
\end{proof}
}
}

{
\lemma{
\label{lemma: regularity and integrability of simple averaged coefficients}
Assume \ref{assumption: positive recurrence}, \ref{assumption: uniform ellipticity} and \hyperlink{H}{$H^{k, 2 + \alpha}$} for $\alpha \in (0, 1)$ and $k \in \N{}_0$. 
Let $b, \sigma, h \in C^{k, 0}$ satisfy for some $C, p > 0$,
\[
\sum_{|\gamma| \leq k} \sup_x \left( | D^\gamma_x b(x, z) | + | D^\gamma_x \sigma(x, z) | + | D^\gamma_x h(x, z) | \right) \leq C ( 1 + |z|^p ). 
\]
Then $\overbar{b}, \overbar{\sigma}, \overbar{a}, \overbar{h} \in C^k_b$.

\begin{proof}
The result follows from Lemma \ref{lemma: regularity and bound on averaged function}.
\end{proof}
}
}

{
\lemma{
\label{lemma: regularity and integrability of poisson averaged coefficients}
Assume \ref{assumption: positive recurrence}, \ref{assumption: uniform ellipticity} and \hyperlink{H}{$H^{1, 2 + \alpha}$} for $\alpha \in (0, 1)$. 
Let $b_\mathrm{I} \in C^{j, \alpha}$, for $j \in \N{}$, be centered with respect to $\mu_\infty$ with the growth condition
\begin{align*}
\left| b_\mathrm{I}(x, z) \right| &\leq C_0 (1 + |z|)^\beta, \\
\sum_{|i| = 1}^j \left| D^i_x b_\mathrm{I}(x, z) \right| &\leq C (1 + |z|^q),
\end{align*}
for \Nicolas{$\beta < -2$}. 
Then $\wt{a} \in C^j_b$ and $\wt{b} \in C^{j - 1}_b$.

\begin{proof}
The result for $\wt{a}$ follows from the fact that the assumptions on $b_\mathrm{I}$ give $\mG_F^{-1}( - b_\mathrm{I} )(\cdot, z) \in C^j_b$ for each $z$ by Theorem \ref{theorem: solution, bounds and regularity of Poisson equation} and the rest then follows from Lemma \ref{lemma: regularity and bound on averaged function}.
For $\wt{b}$ we again have from Theorem \ref{theorem: solution, bounds and regularity of Poisson equation} that $|\nabla_x \mG_F^{-1}( -b_\mathrm{I} ) | \lesssim (1 + |z|^q)$ and $\nabla_x \mG_F^{-1}( -b_\mathrm{I} )(\cdot, z) \in C^{j - 1}$ for each $z$, and therefore we can use Lemma \ref{lemma: regularity and bound on averaged function} to get the desired result. 
\end{proof}
}
}

\subsection{Estimates of the Unnormalized Conditional Distribution}

{
\lemma{
\label{lemma: averaged filter moment bound}
If $\rho^0$ satisfies Eq. \ref{equation: homogenized Zakai}, $\overbar{h}$ is bounded, and $\varphi \in C^2_b(\R{m}; \R{})$,  then for $p \geq 2$, 
\begin{align*}
\E*[\Q{}] \sup_{t \leq T} | \rho^0_t(\varphi) |^p < \infty.
\end{align*}
If $\rho^\epsilon$ satisfies Eq. \ref{equation: Zakai} and $h$ is bounded, then for $p \geq 2$, 
\begin{align*}
\sup_{\epsilon \in (0, 1]} \E*[\Q{}] \sup_{t \leq T} | \rho^\epsilon_t(1) |^p < \infty.
\end{align*}

\begin{proof}
For the first result, since $\varphi$ is bounded we have $\rho^0_t(\varphi) \leq | \varphi |_\infty \rho^0_t(1)$ and therefore we aim to show $\E*[\Q{}] \sup_{t \leq T} | \rho^0_t(1) |^p < \infty$.
Applying $\E*[\Q{}] \sup_{t \leq T} | \cdot |^p$ to the evolution equation for $\rho^0_t(1)$ gives
\begin{align*}
\E*[\Q{}] \sup_{t \leq T} | \rho^0_t(1) |^p 
&\lesssim \E*[\Q{}] \left| \rho^0_0(1) \right|^p 
+ \E*[\Q{}] \sup_{t \leq T} \left| \int_0^t \langle \rho^0_s(\varphi \overbar{h}) ,  d\Ye_s \rangle \right|^p \\
&\lesssim 1 + T^{(p/2) - 1} \int_0^T \E*[\Q{}] | \rho^0_s( \varphi \overbar{h} ) |^p ds,
\end{align*}
where in the last step we applied BDG, H\"older's inequality, and Fubini. 
Using the boundedness of $\varphi \overbar{h}$ and Lemma \ref{lemma: moments of Girsanov}, the integrand is bounded by, 
\begin{align*}
\E*[\Q{}] | \rho^0_s(\varphi \overbar{h}) |^p
&\leq \E[\Q{}]{ \E[\Pe]{ | \varphi \overbar{h}(X^0_s) \wt{D}^0_s  |^p }[\mYe_s]  } \\
&= \E[\Pe]{ \wt{D}^0_s \E[\Pe]{ | \varphi \overbar{h}(X^0_s) \wt{D}^0_s  |^p }[\mYe_s]  } 
\leq \frac{1}{2} \E*[\Pe]{ \left( \wt{D}^0_s \right)^2 } + \frac{1}{2} | \varphi \overbar{h} |_\infty^{2p} \E*[\Pe]{ \left( \wt{D}^0_s \right)^{2p} } < \infty.
\end{align*}
The second result follows in the same manner.
\end{proof}
}
}

\section{Existence, Characterization and Uniqueness of Weak Limits}
\label{section: existence, characterization, and uniqueness}

Let $S(\R{m})$ be the space of finite signed Borel measures on $\R{m}$ with the weak topology induced by $C_b(\R{m}; \R{})$, and $C([0, T]; S(\R{m}))$ the space of continuous paths with values in $S(\R{m})$ endowed with the topology of uniform convergence. 
For each $\epsilon$, we denote the $C([0, T]; S(\R{m}))$-valued random variable $\zeta^\epsilon \equiv \rho^{\epsilon, x} - \rho^0$, the difference of the $x$-marginal and averaged filter.
With this notation, now define the $\epsilon$-parameterized family of Borel probability measures $(P^\epsilon)$ on $C([0, T]; S(\R{m}))$, to be those induced by $(\zeta^\epsilon)$,
\[
P^\epsilon ( \cdot ) = \Q{} \left( \left( \zeta^\epsilon \right)^{-1} ( \cdot ) \right).
\]

We will need to prove a uniform concentration condition\footnote{This is sometimes referred to as compact confinement condition in the literature.} of the collection $(P^\epsilon)$ for the proof of the existence of weak limits of $(\zeta^\epsilon)$.
In the context of our problem, the uniform concentration condition is the following:
{
\defn[Uniform Concentration Condition]{
\label{definition: Jakubowski uniform concentration condition}
$(P^\epsilon)$ is said to satisfy the uniform concentration condition if for each $\delta > 0$, there exists a compact set $K_\delta \subset S(\R{m})$ such that 
\begin{align}
\label{equation: Jakubowski uniform concentration condition}
P^\epsilon (C([0, T]; K_\delta)) \geq 1 - \delta, \quad \forall \epsilon.
\end{align}
}
} 

We now prove a lemma that provides a sufficient condition for the uniform concentration condition.

{
\lem{
\label{lemma: sufficient condition for uniform concentration}
The uniform concentration condition holds if for some $p > 0$ and continuous $M: \R{m} \rightarrow (0, \infty)$ with $\lim_{ |x| \rightarrow \infty} M(x) = \infty$, we have
\begin{align}
\label{equation: assumption for showing CCC}
\sup_{\epsilon \in (0, 1]} \E*[\Q{}] \sup_{t \leq T} \left( |\zeta^\epsilon_t|\left( M \right) \right)^p < \infty.
\end{align}
Here $| \zeta^\epsilon_t |$ is the total variation measure of $\zeta^\epsilon_t$. 
\begin{proof} 
We first show that for $C > 0$, the set 
\begin{align*}
K = \set{\mu \in S(\R{m})}[ |\mu| (M) \leq C ]
\end{align*}
is tight. 
Given $\delta > 0$, choose $R > 0$ large enough such that $\inf_{|x| \geq R} M(x) \geq C / \delta$. 
Then denoting $A_\delta = B(0, R) \subset \R{m}$, the closed ball centered at the origin with radius $R$, we have that for any $\mu \in K$
\begin{align*}
|\mu| (A^c_\delta) 
= |\mu| (1_{|\cdot| > R})
\leq |\mu| \left( \left( \frac{\delta}{C} M \right) 1_{|\cdot| > R} \right)
\leq \frac{\delta}{C} |\mu|(M) 
\leq \frac{\delta}{C} C 
= \delta. 
\end{align*}
This shows that $K$ is tight.
Moreover, since $M$ is bounded from below by some $m > 0$, we have 
\[
\sup_{\mu \in K} | \mu | (1) \leq \frac{1}{m} \sup_{\mu \in K} | \mu | (M) \leq \frac{C}{m}, 
\]
and therefore $K$ is bounded in total variation norm. 
Since $S(\R{m})$ with the weak topology induced by $C_b(\R{m}; \R{})$ is Polish, we have by Prokhorov's theorem that $K$ is relatively compact. 
Further, by Fatou's lemma for weak convergence, $K$ is also closed and therefore compact. 

Given $\delta > 0$, choose $C > 0$ large enough so that 
\begin{align*}
\frac{\sup_{\epsilon \in (0, 1]} \E*[\Q{}] \sup_{t \leq T} \left( |\zeta^\epsilon_t|\left( M \right) \right)^p}{C^p} < \delta.
\end{align*}
Defining our compact set $K_\delta = \set{\mu \in S(\R{m})}[ |\mu| (M) \leq C ]$, we have 
\begin{align*}
\Q{} \left( \zeta^\epsilon \notin C([0, T]; K_\delta) \right) 
\leq \Q{} \left( \sup_{t \leq T} | \zeta^\epsilon_t | (M) > C \right)
\leq \frac{\sup_{\epsilon \in (0, 1]} \E*[\Q{}] \sup_{t \leq T} \left( |\zeta^\epsilon_t|\left( M \right) \right)^p}{C^p} 
\leq \delta.
\end{align*}
\end{proof}
}
}

The next result uses Lemma \ref{lemma: sufficient condition for uniform concentration} to prove that $(P^\epsilon)$ is tight.

{
\lem{
\label{lemma: solutions to signed Zakai equations are tight}
Assume that $f$ and $g$ satisfy the assumptions of \ref{assumption: positive recurrence}, \ref{assumption: uniform ellipticity} and \hyperlink{H}{$H^{2, 2 + \alpha}$} for $\alpha \in (0, 1)$, and that $b, \sigma, h$ are bounded. 
Let $b_\mathrm{I} \in C^{2, \alpha}$, be centered with respect to $\mu_\infty(x)$ for each $x$, and satisfy the growth conditions 
\begin{align*}		
\left| b_\mathrm{I}(x, z) \right| &\leq C (1 + |z|)^\beta, \\
\sum_{|i| = 1}^2 \left| D^i_x b_\mathrm{I}(x, z) \right| &\leq C (1 + |z|^q),
\end{align*}
for some \Nicolas{$\beta < -2$} and $q > 0$.
Assume that $\Q{}_{(\Xe_0, \Ze_0)}$ has finite moments of every order.
Then the $\epsilon$-parameterized family of Borel probability measures $(P^\epsilon)$ 
is tight. 

\begin{proof} 
To prove the statement, we follow criteria provided in \cite[Theorem 3.1, p.276]{Jakubowski:1986}, which gives conditions for a family of Borel probability measures on $D([0, T]; E)$, c\`adl\`ag path space with $E$ a completely regular topological space with metrizable compacts, to be tight. 
$C([0, T]; S(\R{m}))$ is viewed as a subset of $D([0, T]; S(\R{m}))$, and $S(\R{m})$ with the weak topology induced by $C^2_b(\R{m}; \R{})$ is Polish and therefore a completely regular topological space with metrizable compacts.

Specifically, let $\mathbb{F}$ be the natural injection of $C^2_b(\R{m}; \R{})$ into its double dual. 
This collection 
satisfies criteria for \cite[Theorem 3.1, p.276]{Jakubowski:1986}, i.e., it is a collection of continuous functions that separate points in $S(\R{m})$, and is closed under addition (i.e., $f, g \in \mathbb{F}$, then $f + g \in \mathbb{F}$). 
Then to each $f \in \mathbb{F}$ associate a map $\wt{f} \in \wt{\mathbb{F}}$, characterized as follows, 
\begin{align*}
\wt{f} : C([0, T]; S(\R{m})) &\longrightarrow C([0, T]; \R{}) \\
\mu &\longmapsto f \circ \mu.
\end{align*}
The conditions for tightness by \cite[Theorem 3.1, p.276]{Jakubowski:1986} then states that $(P^\epsilon)$ is tight if and only if the following two conditions are satisfied: 
\begin{enumerate}[(i)]
\item For each $\delta > 0$ there is a compact set $K_\delta \subset S(\R{m})$ such that 
\[
P^\epsilon (C([0, T]; K_\delta)) > 1 - \delta, \quad \forall \epsilon
\]
\item The family $(P^\epsilon)$ is $\mathbb{F}$-weakly tight, i.e., for each $f \in \mathbb{F}$ the family $(P^\epsilon \circ (\wt{f}^{-1}))$ of probability measures on $C([0, T]; \R{})$ is tight. 
\end{enumerate}

The proof of (i) will follow from Lemma \ref{lemma: sufficient condition for uniform concentration}, which provides a sufficient condition for the uniform concentration condition. 
To prove the lemma, we define $M(x) = (1 + |x|^2)^{1/2}$, let $p \geq 2$ and check that the following condition holds,
\begin{align}
\tag{\ref{equation: assumption for showing CCC}}
\sup_{\epsilon \in (0, 1]} \E*[\Q{}] \sup_{t \leq T} \left( |\zeta^\epsilon_t|\left( M \right) \right)^p < \infty.
\end{align}
Note that $M(x)$ satisfies the conditions for Lemma \ref{lemma: sufficient condition for uniform concentration} and has bounded first and second order derivatives, 
\begin{align*}
\left| \frac{\partial}{\partial x_i} M(x) \right| 
&= \left| \frac{x_i}{M(x)} \right| 
\lesssim 1, \\
\left| \frac{\partial^2}{\partial x_i \partial x_j} M(x) \right| 
&= \left| -\frac{x_i x_j}{M(x)^3} + \frac{\delta_{ij}}{M(x)} \right|
\lesssim \frac{1}{M(x)} \lesssim 1.
\end{align*}
Directly estimating, we have 
\begin{align*}
\sup_{\epsilon \in (0, 1]} \E*[\Q{}] \sup_{t \leq T} \left( |\zeta^\epsilon_t| (M) \right)^p
= \sup_{\epsilon \in (0, 1]} \E*[\Q{}] \sup_{t \leq T} \left( \rho^{\epsilon, x}_t(M) + \rho^0_t(M) \right)^p 
\lesssim \sup_{\epsilon \in (0, 1]} \E*[\Q{}] \sup_{t \leq T} \left| \rho^{\epsilon, x}_t(M) \right|^p + \E*[\Q{}] \sup_{t \leq T} \left| \rho^0_t(M) \right|^p .
\end{align*}

Dealing with each term separately, we first address $\sup_{\epsilon \in (0, 1]} \E*[\Q{}] \sup_{t \leq T} \left( \rho^{\epsilon, x}_t(M) \right)^p$.
To handle the singular term (the intermediate drift) in the slow process, we perturb $M$ by a corrector term. 
Define the perturbed test function, 
\[
\Me(x, z) = M(x) + \epsilon \chi(x, z),
\]
with $\chi(x, z)$ the solution of the Poisson equation, 
\begin{align*}
\mG_F \chi = - \mG_I M. 
\end{align*}
The Poisson equation is well-posed with the right hand side satisfying the conditions of Theorem \ref{theorem: solution, bounds and regularity of Poisson equation} (recall that $b_\mathrm{I}$ has the correct decay in the $z$ variable), and therefore the regularity and bounds of $\chi$ come from Theorem \ref{theorem: solution, bounds and regularity of Poisson equation}. 
Specifically, 
\begin{align*}
|\chi(x, z)| &\lesssim 1, \\
\sum_{|i| = 1}^2 |D^i_x \chi(x, z)| &\lesssim  (1 + |z|^{q'}), 
\end{align*}
for some $q' > 0$. 
Using the identity $\rho^\epsilon_t(\Me) = \rho^{\epsilon, x}_t(M) + \epsilon \rho^\epsilon_t(\chi)$, we have the representation 
\begin{align*}
\rho^{\epsilon, x}_t(M)
= -\epsilon \rho^\epsilon_t(\chi) 
&+ \rho^\epsilon_0(\Me) + \int_0^t \rho^\epsilon_s(\mG_S M) ds + \int_0^t \rho^\epsilon_s (\mG_\mathrm{I} \chi) ds + \epsilon \int_0^t \rho^\epsilon_s ( \mG_S \chi) ds \\
&+ \int_0^t \langle \rho^\epsilon_s(M h + \alpha \sigma^* \nabla_x M ), d\Ye_s \rangle + \epsilon \int_0^t \langle \rho^\epsilon_s(\chi h + \alpha \sigma^* \nabla_x \chi ), d\Ye_s \rangle. 
\end{align*}
And therefore, 
\EquationAligned{
\label{equation: starting point for uniform concentration proof}
\E*[\Q{}] \sup_{t \leq T} \left| \rho^{\epsilon, x}_t(M) \right|^p
&\lesssim \epsilon^p \E*[\Q{}]  \sup_{t \leq T} \left| \rho^\epsilon_t(\chi) \right|^p
+ \E*[\Q{}] \left| \rho^\epsilon_0(\Me) \right|^p \\
&+ \E*[\Q{}] \int_0^T \left| \rho^\epsilon_s(\mG_S M) \right|^p ds  
+ \E*[\Q{}] \int_0^T \left| \rho^\epsilon_s (\mG_\mathrm{I} \chi) \right|^p ds 
+ \epsilon^p \E*[\Q{}] \int_0^T \left| \rho^\epsilon_s ( \mG_S \chi) \right|^p ds \\
&+ \E*[\Q{}] \sup_{t \leq T} \left| \int_0^t \langle \rho^\epsilon_s(M h), d\Ye_s \rangle \right|^p 
+ \E*[\Q{}] \sup_{t \leq T} \left| \int_0^t \langle \rho^\epsilon_s( \alpha \sigma^* \nabla_x M ), d\Ye_s \rangle \right|^p \\
&+ \epsilon^p \E*[\Q{}] \sup_{t \leq T} \left| \int_0^t \langle \rho^\epsilon_s(\chi h), d\Ye_s \rangle \right|^p 
+ \epsilon^p \E*[\Q{}] \sup_{t \leq T} \left| \int_0^t \langle \rho^\epsilon_s( \alpha \sigma^* \nabla_x \chi ), d\Ye_s \rangle \right|^p.
}

By the boundedness of $\chi$ and application of Lemma \ref{lemma: averaged filter moment bound}, the first term of Eq. \ref{equation: starting point for uniform concentration proof} is
\begin{align*}
\sup_{\epsilon \in (0, 1]} \E*[\Q{}] \sup_{t \leq T} \left| \rho^\epsilon_t(\chi) \right|^p
\lesssim \sup_{\epsilon \in (0, 1]} \E*[\Q{}] \sup_{t \leq T} \rho^\epsilon_t( \left| \chi \right| )^p 
\lesssim \sup_{\epsilon \in (0, 1]} \E*[\Q{}] \sup_{t \leq T} \rho^\epsilon_t( 1 )^p 
< \infty.
\end{align*}
For the second term of Eq. \ref{equation: starting point for uniform concentration proof}, we have
\begin{align*}
\E*[\Q{}] | \rho^\epsilon_0(M^\epsilon) |^p
&= \E*[\Q{}] \left| \E[\Pe]{M^\epsilon(\Xe_0, \Ze_0) \wt{D}^\epsilon_0}[\mYe_0] \right|^p
= \E*[\Q{}] { \left| \E[\Q{}] { M^\epsilon(\Xe_0, \Ze_0) } \right|^p } \\
&\leq \int | M^\epsilon(x, z) |^p \Q{}_{(\Xe_0, \Ze_0)}(dx, dz)
\lesssim \int M(x)^p + \epsilon^p | \chi(x, z) |^p \Q{}_{(\Xe_0, \Ze_0)}(dx, dz), 
\end{align*}
and therefore 
\begin{align*}
\sup_{\epsilon \in (0, 1]} \E*[\Q{}] | \rho^\epsilon_0(M^\epsilon) |^p < \infty.
\end{align*}
Where we used the boundedness of $\chi$ and the fact that we have finite moments of every order for $\Q{}_{(\Xe_0, \Ze_0)}$.
From the boundedness of $b, \sigma, D^\gamma_x M$, for $|\gamma| \in \set{1, 2}$, Lemma \ref{lemma: moments of Girsanov} and \ref{lemma: averaged filter moment bound}, we have for the third term of Eq. \ref{equation: starting point for uniform concentration proof}
\begin{align*}
\sup_{\epsilon \in (0, 1]} \E*[\Q{}] \int_0^T \left| \rho^\epsilon_s(\mG_S M) \right|^p ds  
&\leq | \mG_S M |_\infty^p \sup_{\epsilon \in (0, 1]} \E*[\Q{}] \sup_{t \leq T} \rho^\epsilon_s(1)^p ds
< \infty.
\end{align*}
For the fourth term of Eq. \ref{equation: starting point for uniform concentration proof}
\begin{align*}
\E*[\Q{}] \int_0^T \left| \rho^\epsilon_s (\mG_\mathrm{I} \chi) \right|^p ds 
&\lesssim \int_0^T \E[\Pe]{ (\wDe_s)^2 }^{1/2} \E[\Pe]{ \left| \mG_\mathrm{I} \chi (\Xe_s, \Ze_s) \wDe_s \right|^{2p} }^{1/2}  ds \\
&\leq \int_0^T \E[\Pe]{ (\wDe_s)^2 }^{1/2} \E[\Pe]{ \left( \wDe_s \right)^{4p} }^{1/4} \E[\Pe]{ \left| \mG_\mathrm{I} \chi (\Xe_s, \Ze_s) \right|^{4p} }^{1/4} ds.
\end{align*}
And we have 
\begin{align}
\label{equation: existence, fourth term}
\sup_{\epsilon \in (0, 1]} \E[\Pe]{ \left| \mG_\mathrm{I} \chi (\Xe_s, \Ze_s) \right|^{4p} }
\lesssim \sup_{\epsilon \in (0, 1]} \E[\Pe]{ 1 + |\Ze_s|^q } 
\lesssim 1 + \sup_{\epsilon \in (0, 1]} \E[\Q{}]{ |\Ze_0|^q } 
< \infty, 
\end{align}
for some $q > 0$, and therefore 
\begin{align*}
\sup_{\epsilon \in (0, 1]} \E*[\Q{}] \int_0^T \left| \rho^\epsilon_s (\mG_\mathrm{I} \chi) \right|^p ds 
\lesssim \sup_{\epsilon \in (0, 1]} \int_0^T \E[\Pe]{ \left| \mG_\mathrm{I} \chi (\Xe_s, \Ze_s) \right|^{4p} }^{1/4} ds
< \infty.
\end{align*}
By the same arguments as for the fourth term, we have that the fifth term of Eq. \ref{equation: starting point for uniform concentration proof} is bounded uniformly in $\epsilon$, 
\begin{align*}
\sup_{\epsilon \in (0, 1]} \E*[\Q{}] \int_0^T \left| \rho^\epsilon_s ( \mG_S \chi) \right|^p ds < \infty.
\end{align*}
The first stochastic integral, sixth term of Eq. \ref{equation: starting point for uniform concentration proof}, is handled with application of BDG, H\"older's inequality, and Fubini on the first line and then change of measure, Cauchy-Schwarz, and H\"older's inequality to give
\begin{align*}
\sup_{\epsilon \in (0, 1]} \E*[\Q{}] \sup_{t \leq T} \left| \int_0^t \langle \rho^\epsilon_s(M h), d\Ye_s \rangle \right|^p 
&\lesssim \sup_{\epsilon \in (0, 1]} \int_0^T \E*[\Q{}] \left| \rho^\epsilon_s(M h) \right|^p ds \\
&\lesssim \sup_{\epsilon \in (0, 1]} \int_0^T \E[\Pe]{ (\wDe_s)^2 }^{1/2} \E[\Pe]{ \left| \rho^\epsilon_s(M h) \right|^{2p} }^{1/2} ds.
\end{align*}
And by further application of Cauchy-Schwarz, H\"older's inequality, boundedness of $h$, Lemma \ref{lemma: moments on slow process} for some $q > 0$, and finite moments of all orders for $\Q{}_{(\Xe_0, \Ze_0)}$, we get
\begin{align*}
\sup_{\epsilon \in (0, 1]} \E[\Pe]{ \left| \rho^\epsilon_s(M h) \right|^{2p} }^{1/2}
&\leq \sup_{\epsilon \in (0, 1]} \E[\Pe]{ (\wDe_s)^{4p} }^{1/4 } |h|_\infty^p \E[\Pe]{ \left| M(\Xe_s) \right|^{4p} }^{1/4} \\
&\lesssim \sup_{\epsilon \in (0, 1]} \E[\Pe]{ (1 + |\Xe_s|^2)^{2p} }^{1/4} 
\lesssim \sup_{\epsilon \in (0, 1]} \left( 1 + \E[\Pe] { |\Xe_s|^{4p} } \right)^{1/4} \\
&\lesssim \sup_{\epsilon \in (0, 1]} \left( 1 + \E[\Q{}] { |\Xe_0|^{4p} + (1 + \epsilon^{4p}) (1 + |\Ze_0|^q) } \right)^{1/4}
< \infty.
\end{align*}
Because $ \alpha \sigma^* \nabla_x M $ is bounded and by Lemma \ref{lemma: averaged filter moment bound}, the seventh term of Eq. \ref{equation: starting point for uniform concentration proof} is bounded uniformly in $\epsilon$, 
\begin{align*}
\sup_{\epsilon \in (0, 1]} \E*[\Q{}] \sup_{t \leq T} \left| \int_0^t \langle \rho^\epsilon_s( \alpha \sigma^* \nabla_x M ), d\Ye_s \rangle \right|^p 
\lesssim \int_0^T \sup_{\epsilon \in (0, 1]} \E*[\Q{}] \left| \rho^\epsilon_s( \alpha \sigma^* \nabla_x M ) \right|^p ds
< \infty.
\end{align*}
Similarly using the boundedness of $\chi h$, the eighth term of Eq. \ref{equation: starting point for uniform concentration proof} is bounded uniformly in $\epsilon$, 
\begin{align*}
\sup_{\epsilon \in (0, 1]} \E*[\Q{}] \sup_{t \leq T} \left| \int_0^t \langle \rho^\epsilon_s(\chi h), d\Ye_s \rangle \right|^p 
< \infty.
\end{align*}
Making use of the boundedness of $\sigma$, the polynomial growth of $\nabla_x \chi$ in $z$, and the finite moments of all orders for $\Q{}_{(\Xe_0, \Ze_0)}$, the last term of Eq. \ref{equation: starting point for uniform concentration proof} is
\begin{multline*}
\sup_{\epsilon \in (0, 1]} \E*[\Q{}] \sup_{t \leq T} \left| \int_0^t \langle \rho^\epsilon_s( \alpha \sigma^* \nabla_x \chi ), d\Ye_s \rangle \right|^p 
\lesssim \int_0^T \sup_{\epsilon \in (0, 1]} \E*[\Q{}] \left| \rho^\epsilon_s( \alpha \sigma^* \nabla_x \chi ) \right|^p ds \\
\lesssim \int_0^T \sup_{\epsilon \in (0, 1]} \E[\Pe] { \left| \alpha \sigma^* \nabla_x \chi (\Xe_s, \Ze_s) \right|^{4p}  }^{1 / 4} ds 
\lesssim \int_0^T \sup_{\epsilon \in (0, 1]} \E[\Pe] { \left| \nabla_x \chi (\Xe_s, \Ze_s) \right|^{4p}  }^{1 / 4} ds \\
\lesssim \int_0^T \sup_{\epsilon \in (0, 1]} \left( 1 + \E[\Pe] { \left| \Ze_s \right|^q } \right)^{1 / 4} ds 
\lesssim \int_0^T \sup_{\epsilon \in (0, 1]} \left( 1 + \E[\Q{}] { \left| \Ze_0 \right|^q } \right)^{1 / 4} ds
< \infty.
\end{multline*}
This completes the calculation that $\sup_\epsilon \E*[\Q{}] \sup_{t \leq T} \left| \rho^{\epsilon, x}_t(M) \right|^p < \infty$. 
We now show that $\E*[\Q{}] \sup_{t \leq T} \left| \rho^0_t(M) \right|^p < \infty$. 
We have, 
\EquationAligned{
\label{equation: existence averaged term}
\E*[\Q{}] \sup_{t \leq T} \left| \rho^0_t(M) \right|^p
\lesssim \E*[\Q{}] \left| \rho^0_0(M) \right|^p
&+ \int_0^T \E*[\Q{}] \left| \rho^0_s(\mG^\dagger M) \right|^p ds \\
&+ \E*[\Q{}] \sup_{t \leq T} \left| \int_0^t \langle \rho^0_s(M \overbar{h} ), d\Ye_s \rangle \right|^p
+ \E*[\Q{}] \sup_{t \leq T} \left| \int_0^t \langle \rho^0_s( \alpha \overbar{\sigma}^* \nabla_x M ), d\Ye_s \rangle \right|^p.
}
By Lemma \ref{lemma: regularity and integrability of simple averaged coefficients} and \ref{lemma: regularity and integrability of poisson averaged coefficients}, we have that $\overbar{b}, \overbar{\sigma}, \overbar{h}, \wt{a}$ and $\wt{b}$ are all bounded functions.
Therefore, similar arguments as for $\rho^\epsilon$ show that the right side of Eq. \ref{equation: existence averaged term} is finite. 

We now prove (ii), the $\mathbb{F}$-weak tightness condition. 
Let $\varphi \in C^2_b(\R{m}; \R{})$. 
There are two conditions that must be checked for each fixed $\varphi \in \mathbb{F}$. 
The first one is a boundedness condition, 
\begin{align*}
 \lim_{N \rightarrow \infty} \sup_\epsilon \Q{} (| \zeta^\epsilon_0(\varphi) | \geq N) = 0,
\end{align*}
which is trivially satisfied since $\zeta^\epsilon_0(\varphi) = \E[\Q{}]{\varphi(\Xe_0)} - \E[\Q{}]{\varphi(X^0_0)} = 0$.
The second one is an equicontinuity condition\textendash for each $\delta > 0$ and $\eta > 0$ there are $\Delta > 0$ and $j < \infty$ such that 
\begin{align}
\label{equation: necessary equicontinuity condition for weak tightness}
\Q{} \left( \sup_{|t - s| \leq \Delta} |\zeta^{\epsilon_i}_t(\varphi) - \zeta^{\epsilon_i}_s(\varphi) | \geq \delta \right) \leq \eta, \quad \forall i \geq j.
\end{align}
A sufficient condition for Eq. \ref{equation: necessary equicontinuity condition for weak tightness} is the following\textendash there are $\mu, \beta, \gamma > 0$ and $K < \infty$ such that
\begin{align}
\label{equation: sufficient equicontinuity condition for weak tightness}
\E*[\Q{}] | \zeta^{\epsilon_i}_t(\varphi) - \zeta^{\epsilon_i}_s(\varphi) |^\mu \leq K |t - s|^{1 + \beta} + \epsilon^\gamma_i, \quad \forall i.
\end{align}
We now show that Eq. \ref{equation: sufficient equicontinuity condition for weak tightness} is true.
Let $\mu = 4$, $\epsilon > 0$ and $\varphi^\epsilon = \varphi + \epsilon \chi$, where $\chi$ solves the Poisson equation, 
\begin{align*}
\mG_F \chi = - \mG_\mathrm{I} \varphi. 
\end{align*}
From $\zeta^\epsilon_t(\varphi) = \rho^{\epsilon, x}_t(\varphi) - \rho^0_t(\varphi)$ and $\rho^\epsilon_t(\varphi^\epsilon) = \rho^{\epsilon, x}_t(\varphi) + \rho^\epsilon_t(\epsilon \chi)$ we have 
\begin{align}
\label{equation: expansion of sufficient equicontinuity condition in F-weak tightness}
| \zeta^\epsilon_t(\varphi) - \zeta^\epsilon_s(\varphi) |^4 
\lesssim | \rho^\epsilon_t(\varphi^\epsilon) - \rho^\epsilon_s(\varphi^\epsilon) |^4 + | \rho^0_t(\varphi) - \rho^0_s(\varphi) |^4 + | \rho^\epsilon_t (\epsilon \chi) - \rho^\epsilon_s (\epsilon \chi) |^4. 
\end{align}
Working on the first term in this inequality, 
\begin{multline}
\label{equation: expansion of condition in equicontinuity proof}
\E*[\Q{}] | \rho^\epsilon_t(\varphi^\epsilon) - \rho^\epsilon_s(\varphi^\epsilon) |^4
\lesssim \E*[\Q{}] \left| \int_s^t \rho^\epsilon_u ( \mG_S \varphi + \mG_\mathrm{I} \chi + \epsilon \mG_S \chi)  du \right|^4 + \E*[\Q{}] \left| \int_s^t \langle \rho^\epsilon_u (\varphi^\epsilon h + \alpha \sigma^* \nabla_x \varphi^\epsilon ), d\Ye_u \rangle \right|^4 \\
\lesssim (t - s)^3 \int_s^t \E*[\Q{}] \left| \rho^\epsilon_u ( \mG_S \varphi + \mG_\mathrm{I} \chi + \epsilon \mG_S \chi)  \right|^4 du 
+ (t - s) \int_s^t \E*[\Q{}] \left| \rho^\epsilon_u (\varphi^\epsilon h + \alpha \sigma^* \nabla_x \varphi^\epsilon ) \right|^4 du. 
\end{multline}
First term of Eq. \ref{equation: expansion of condition in equicontinuity proof}, 
\begin{align*}
\int_s^t \E*[\Q{}] \left| \rho^\epsilon_u ( \mG_S \varphi + \mG_\mathrm{I} \chi + \epsilon \mG_S \chi) \right|^4 du 
\lesssim \int_s^t \E*[\Q{}] \left| \rho^\epsilon_u ( \mG_S \varphi) \right|^4 + \E*[\Q{}] \left| \rho^\epsilon_u ( \mG_\mathrm{I} \chi ) \right|^4 + \E*[\Q{}] \left| \rho^\epsilon_u ( \epsilon \mG_S \chi) \right|^4 du.
\end{align*}
By the boundedness of $b, \sigma, D^k_x \varphi$ for $|k| \leq 2$ and Lemma \ref{lemma: moments of Girsanov}, we have the first term bounded by
\begin{align*}
\int_s^t \E*[\Q{}] \left| \rho^\epsilon_u ( \mG_S \varphi) \right|^4 du 
\lesssim (t - s).
\end{align*}
By the boundedness of $b, \sigma$, the polynomial growth of $b_\mathrm{I}, D^k_x \chi$ in $z$, for $|k| \leq 2$ and Lemma \ref{lemma: moments of Girsanov}, we have by the same arguments as Eq. \ref{equation: existence, fourth term}
\begin{align*}
\int_s^t \E*[\Q{}] \left| \rho^\epsilon_u ( \mG_\mathrm{I} \chi ) \right|^4 du 
+ \int_s^t \E*[\Q{}] \left| \rho^\epsilon_u ( \epsilon \mG_S \chi) \right|^4 du
\lesssim (t - s)
+ \epsilon^4 (t - s)
= (1 + \epsilon^4) (t - s). 
\end{align*}
And therefore we have, 
\begin{align*}
\int_s^t \E*[\Q{}] \left| \rho^\epsilon_u ( \mG_S \varphi + \mG_\mathrm{I} \chi + \epsilon \mG_S \chi) \right|^4 du 
\lesssim (1 + \epsilon^4) (t - s) .
\end{align*}
The second term of Eq. \ref{equation: expansion of condition in equicontinuity proof} is bounded as follows, 
\begin{align*}
\int_s^t \E*[\Q{}] \left| \rho^\epsilon_u (\varphi^\epsilon h + \alpha \sigma^* \nabla_x \varphi^\epsilon ) \right|^4 du
&\lesssim \int_s^t \E*[\Q{}] \left| \rho^\epsilon_u (\varphi h) \right|^4 + \E*[\Q{}] \left| \rho^\epsilon_u ( \epsilon \chi h ) \right|^4 du \\
&+ \int_s^t \E*[\Q{}] \left| \rho^\epsilon_u ( \alpha \sigma^* \nabla_x \varphi ) \right|^4 + \E*[\Q{}] \left| \rho^\epsilon_u ( \epsilon \alpha \sigma^* \nabla_x \chi ) \right|^4  du \\
&\lesssim (1 + \epsilon^4) (t - s) + \E*[\Q{}] \left| \rho^\epsilon_u ( \epsilon \alpha \sigma^* \nabla_x \chi ) \right|^4  du, 
\end{align*}
by boundedness of $h, \sigma, \chi, D^k_x \varphi$ for $|k| \leq 2$, and Lemma \ref{lemma: moments of Girsanov}.
The last term is bounded by $\lesssim \epsilon^4 (t - s)$ due to boundedness of $\sigma$ and polynomial growth of $\nabla_x \chi$ in $z$ and finite moments of all orders for $\Q{}_{(\Xe_0, \Ze_0)}$, and therefore 
\begin{align*}
\int_s^t \E*[\Q{}] \left| \rho^\epsilon_u (\varphi^\epsilon h + \alpha \sigma^* \nabla_x \varphi^\epsilon ) \right|^4 du
\lesssim (1 + \epsilon^4) (t - s) .
\end{align*}
The expectation of the second term in Eq. \ref{equation: expansion of sufficient equicontinuity condition in F-weak tightness} is, 
\begin{align*}
\E*[\Q{}] | \rho^0_t(\varphi) - \rho^0_s(\varphi) |^4 
&\lesssim \E*[\Q{}] \left| \int_s^t \rho^0_u(\mG^\dagger \varphi) du \right|^4 
+ \E*[\Q{}] \left| \int_s^t \langle \rho^0_u( \varphi \overbar{h} + \alpha \overbar{\sigma}^* \nabla_x \varphi ), d\Ye_u \rangle \right|^4 \\
&\lesssim (t - s)^3 \int_s^t \E*[\Q{}] \left| \rho^0_u(\mG^\dagger \varphi) \right|^4 du + (t - s) \int_s^t \E*[\Q{}] \left| \rho^0_u ( \varphi \overbar{h} + \alpha \overbar{\sigma}^* \nabla_x \varphi ) \right|^4 du \\
&\lesssim (t - s)^4 + (t - s)^2, 
\end{align*}
where we use the boundedness of $\overbar{h}, \overbar{\sigma}$ and the coefficients of $\mG^\dagger$, which is a result of Lemma \ref{lemma: regularity and integrability of simple averaged coefficients} and \ref{lemma: regularity and integrability of poisson averaged coefficients}, and the boundedness of $D^k_x \varphi$ for $k \leq 2$.

Using the fact that $\chi$ is bounded and Lemma \ref{lemma: averaged filter moment bound}, the last term of Eq. \ref{equation: expansion of sufficient equicontinuity condition in F-weak tightness} is bounded by
\begin{align*}
\E*[\Q{}] | \rho^\epsilon_t (\epsilon \chi) - \rho^\epsilon_s (\epsilon \chi) |^4
\lesssim \epsilon^4 \E[\Q{}] { \rho^\epsilon_t (1)^4 + \rho^\epsilon_s (1)^4 }
\lesssim \epsilon^4.
\end{align*}
Lastly, since we are interested in the case $(t -s) < 1$, we satisfy Eq. \ref{equation: sufficient equicontinuity condition for weak tightness} with $\mu = 4$, $\beta = 1$ and $\gamma = 4$. 

\end{proof}
}
}

Now that $(P^\epsilon)$ has been shown to be tight, we need to show that this collection is weakly relatively compact and therefore any subsequence 
will converge to a weak limit. 
The fact that $(P^\epsilon)$ is weakly relatively compact is given by the following corollary.

{
\cor{
\label{corollary: solutions to signed Zakai equations are relatively compact}
Every subsequence of the $\epsilon$-parameterized family of probability measures $(P^\epsilon)$ induced on path space $C([0, T]; S(\R{m}))$ by $\zeta^\epsilon$, has a weak limit.
\begin{proof}
From Lemma \ref{lemma: solutions to signed Zakai equations are tight}, $(P^\epsilon)$ is tight. 
Because $C([0, T]; S(\R{m}))$ is Hausdorff, since it is metrizable by \cite[Proposition 1.6iii, p.267]{Jakubowski:1986} because $S(\R{m})$ is Polish, this implies that $(P^\epsilon)$ is relatively compact (see for instance \cite[Theorem 2.2.1, p.56]{Kallianpur:1995}).
Therefore each sequence of $(P^\epsilon)$ has a convergent subsequence. 
\end{proof}
}
}

Given a subsequence of $(\epsilon)$, we now let $\zeta$ be the limit point and characterize this limit point in the next lemma.

{
\lem{
\label{lemma: characterization of weak limits}
Assume that $f, g$ satisfy the assumptions of \ref{assumption: positive recurrence}, \ref{assumption: uniform ellipticity} and \hyperlink{H}{$H^{2, 2 + \alpha}$}.
Let $b, b_\mathrm{I}, a \in C^{2, \alpha}$ and satisfy the growth conditions 
\begin{align*}		
\left| b(x, z) \right| + \left| b_\mathrm{I}(x, z) \right| + \left| a(x, z) \right| &\leq C (1 + |z|)^\beta, \\
\sum_{|k| = 1}^2 \left| D^k_x b(x, z) \right| + \left| D^k_x b_\mathrm{I}(x, z) \right| + \left| D^k_x a(x, z) \right| &\leq C (1 + |z|^q),
\end{align*}
for some \Nicolas{$\beta < -2$} and $q > 0$.
Let $b_\mathrm{I}$ be centered with respect to $\mu_\infty(x)$ for each $x$.
Assume $h$ is bounded and globally Lipschitz in $(x, z)$.
Let $\sigma$ be globally Lipschitz in $z$.\ifbool{ShowFootnotes}{Note that $\sigma$ is globally Lip. in $x$ as well from the fact that $a$ has bounded first derivatives in $x$}{}
And assume that $\Q{}_{(\Xe_0, \Ze_0)}$ has finite moments of every order.
Then any limit point $\zeta$ of $(\zeta^\epsilon)$ satisfies the equation, 
\begin{align*}
\zeta_t(\varphi) = \int_0^t \zeta_s ( \mG^\dagger \varphi ) ds + \int_0^t \langle \zeta_s ( \varphi \overbar{h} + \alpha \overbar{\sigma}^* \nabla_x \varphi ),  dY_s \rangle, 
\qquad \zeta_0(\varphi) = 0, 
\qquad \text{$\Q{}$-a.s. uniformly in $t \in [0, T]$}.
\end{align*}
\begin{proof} 
With an abuse of notation, let $\epsilon$ be an element of the subsequence $(\epsilon)$, assume $\varphi \in C^2_b(\R{m}; \R{})$, and consider the perturbed test function, 
\begin{align*}
\varphi^\epsilon(x, z) = \varphi(x) + \epsilon \chi(x, z) + \epsilon^2 \wt{\chi}(x, z), 
\end{align*}
where $\chi$ and $\wt{\chi}$ solve the Poisson equations, 
\begin{align*}
\mG_F \chi &= - \mG_\mathrm{I} \varphi, \\
\mG_F \wt{\chi} &= - \left( \mG_S - \overbar{\mG_S} \right) \varphi - \left( \mG_\mathrm{I} \chi - \overbar{\mG_\mathrm{I} \chi} \right). 
\end{align*}
From Theorem \ref{theorem: solution, bounds and regularity of Poisson equation}, we have  
\begin{align*}
|\chi(x, z)| + |\wt{\chi}(x, z)| &\lesssim 1, \\
\sum_{|i| = 1}^2 |D^i_x \chi(x, z)| + | D^i_x \wt{\chi}(x, z)| &\lesssim  (1 + |z|^{\Nicolas{q'}}), 
\end{align*}
for some \Nicolas{$q' > 0$}. 
Because $\rho^{\epsilon, x}_t(\varphi) = \rho^\epsilon_t(\varphi^\epsilon) - \rho^\epsilon_t(\epsilon \chi) - \rho^\epsilon_t(\epsilon^2 \wt{\chi})$, we have
\begin{align}
\label{equation: evolutionary equation for signed measure}
\zeta^{\epsilon}_t(\varphi) = 
-\rho^\epsilon_t(\epsilon \chi) - \rho^\epsilon_t(\epsilon^2 \wt{\chi}) 
&+ \rho^\epsilon_0(\varphi^\epsilon) - \rho^0_0(\varphi)
+ \int_0^t \rho^{\epsilon}_s ( \mG^{\epsilon} \varphi^{\epsilon} ) ds - \int_0^t \rho^0_s ( \mG^\dagger \varphi ) ds \nonumber \\
&+ \int_0^t \langle \rho^{\epsilon}_s ( \varphi^\epsilon h + \alpha \sigma^* \nabla_x \varphi^\epsilon ),  d\Ye_s \rangle - \int_0^t \langle \rho^0_s ( \varphi \overbar{h} + \alpha \overbar{\sigma}^* \nabla_x \varphi ),  d\Ye_s \rangle.
\end{align}
When expanded, the Lebesgue integral for $\rho^{\epsilon}_s ( \mG^{\epsilon} \varphi^{\epsilon} )$ becomes, 
\begin{align*}
\int_0^t \rho^{\epsilon}_s ( \mG^{\epsilon} \varphi^{\epsilon} ) ds 
&= \int_0^t \rho^{\epsilon, x}_s ( \overbar{\mG_S} \varphi ) ds 
+ \int_0^t \rho^{\epsilon, x}_s ( \overbar{\mG_\mathrm{I} \chi} ) ds \\
&+ \epsilon \int_0^t \rho^\epsilon_s ( \mG_S \chi ) ds 
+ \epsilon \int_0^t \rho^\epsilon_s ( \mG_\mathrm{I} \wt{\chi} ) ds 
+ \epsilon^2 \int_0^t \rho^\epsilon_s ( \mG_S \wt{\chi} ) ds. 
\end{align*}
The term $\rho^\epsilon_0(\varphi^\epsilon) - \rho^0_0(\varphi)$ is, 
\begin{align*}
\rho^\epsilon_0(\varphi^\epsilon) - \rho^0_0(\varphi) 
&= \rho^{\epsilon, x}_0(\varphi) + \rho^\epsilon_0( \epsilon \chi ) + \rho^\epsilon_0 ( \epsilon^2 \wt{\chi} ) - \rho^0_0(\varphi)  \\
&= \rho^\epsilon_0( \epsilon \chi ) + \rho^\epsilon_0 ( \epsilon^2 \wt{\chi} ).
\end{align*}
And we group all terms of first order in $\epsilon$ involving $\chi$ into $\mathcal{O}_\chi(\epsilon)$, 
\begin{align*}
\mathcal{O}_\chi(\epsilon) = 
-\rho^\epsilon_t(\epsilon \chi) 
+ \rho^\epsilon_0( \epsilon \chi ) 
+ \epsilon \int_0^t \rho^\epsilon_s ( \mG_S \chi ) ds 
+ \epsilon \int_0^t \langle \rho^{\epsilon}_s ( \chi  h + \alpha \sigma^* \nabla_x \chi ),  d\Ye_s \rangle.
\end{align*}
Similarly, let the terms of first and second order in $\epsilon$ involving $\wt{\chi}$ be grouped into $\mathcal{O}_{\wt{\chi}}(\epsilon)$, 
\begin{align*}
\mathcal{O}_{\wt{\chi}}(\epsilon) = 
- \rho^\epsilon_t(\epsilon^2 \wt{\chi}) 
+ \rho^\epsilon_0 ( \epsilon^2 \wt{\chi} ) 
+ \epsilon \int_0^t \rho^\epsilon_s ( \mG_\mathrm{I} \wt{\chi} ) ds
+ \epsilon^2 \int_0^t \rho^\epsilon_s ( \mG_S \wt{\chi} ) ds
+ \epsilon^2 \int_0^t \langle \rho^{\epsilon}_s ( \wt{\chi} h + \alpha \sigma^* \nabla_x \wt{\chi} ),  d\Ye_s \rangle.
\end{align*}
Eq. \ref{equation: evolutionary equation for signed measure} now becomes, 
\begin{align}
\label{equation: modified evolutionary equation for signed measure}
\zeta^{\epsilon}_t(\varphi) 
= \int_0^t \zeta^\epsilon_s ( \overbar{\mG_S} \varphi ) ds
&+ \int_0^t \rho^{\epsilon, x}_s ( \overbar{\mG_\mathrm{I} \chi} ) ds \nonumber 
+ \mathcal{O}_\chi(\epsilon)
+ \mathcal{O}_{\wt{\chi}}(\epsilon) \\
&+ \int_0^t \langle \rho^{\epsilon}_s ( \varphi h + \alpha \sigma^* \nabla_x \varphi ),  d\Ye_s \rangle - \int_0^t \langle \rho^0_s ( \varphi \overbar{h} + \alpha \overbar{\sigma}^* \nabla_x \varphi ),  d\Ye_s \rangle.
\end{align}
Next, consider the equivalence of the following terms
\begin{align*}
\int_0^t \rho^{\epsilon, x}_s ( \overbar{\mG_\mathrm{I} \chi} ) ds
= \int_0^t \rho^{\epsilon, x}_s ( \wt{\mG} \varphi ) ds.
\end{align*}
This follows since $\nabla^{\otimes 2}_x \varphi = \nabla_x \nabla_x \varphi$ is symmetric, and therefore, 
\begin{align*}
\langle \nabla^{\otimes 2}_x \varphi \mG_F^{-1}(-b_\mathrm{I}), b_\mathrm{I} \rangle 
&= \sum_{i, j=1}^m \frac{\partial^2 \varphi}{\partial x_i \partial x_j} \mG_F^{-1}(-b_\mathrm{I})_i b_{\mathrm{I}, j } \\
&= \sum_{i, j=1}^m \frac{\partial^2 \varphi}{\partial x_i \partial x_j} \left( \frac{1}{2} (\mG_F^{-1}(-b_\mathrm{I}) \otimes b_\mathrm{I})_{ij} + \frac{1}{2} (b_\mathrm{I} \otimes \mG_F^{-1}(-b_\mathrm{I}) )_{ij} \right),
\end{align*}
which leads to the following, 
\begin{align*}
\overbar{\mG_\mathrm{I} \chi} (x) 
&= \int_{\R{n}} \langle \nabla_x \chi, b_\mathrm{I} \rangle(x, z) \mu_\infty(dz; x)
= \int_{\R{n}} \langle \nabla^{\otimes 2}_x \varphi \mG_F^{-1}(-b_\mathrm{I}) + \left( \nabla_x \mG_F^{-1}(-b_\mathrm{I}) \right)^* \nabla_x \varphi , b_\mathrm{I} \rangle(x, z) \mu_\infty(dz; x) \\
&= \langle \nabla_x \varphi(x), \int_{\R{n}} \left( \nabla_x \mG_F^{-1}(-b_\mathrm{I}) \right) b_\mathrm{I} (x, z) \mu_\infty(dz; x) \rangle
+ \int_{\R{n}} \langle \nabla^{\otimes 2}_x \varphi \mG_F^{-1}(-b_\mathrm{I}), b_\mathrm{I} \rangle(x, z) \mu_\infty(dz; x) \\
&= \wt{\mG} \varphi(x).
\end{align*}
Using this equivalence, adding and subtracting the term, $\int_0^t \langle \rho^{\epsilon, x}_s ( \varphi \overbar{h} + \alpha \overbar{\sigma}^* \nabla_x \varphi ),  d\Ye_s \rangle$, then Eq. \ref{equation: modified evolutionary equation for signed measure} becomes,
\begin{align*}
\zeta^{\epsilon}_t(\varphi) 
=\int_0^t \zeta^\epsilon_s ( \mG^\dagger \varphi ) ds 
&+ \int_0^t \langle \zeta^\epsilon_s ( \varphi \overbar{h} + \alpha \overbar{\sigma}^* \nabla_x \varphi ),  d\Ye_s \rangle 
+ \mathcal{O}_\chi(\epsilon)
+ \mathcal{O}_{\wt{\chi}}(\epsilon) \\
&+ \int_0^t \langle \rho^{\epsilon}_s ( \varphi h + \alpha \sigma^* \nabla_x \varphi ), d\Ye_s \rangle - \int_0^t \langle \rho^{\epsilon, x}_s ( \varphi \overbar{h} + \alpha \overbar{\sigma}^* \nabla_x \varphi ), d\Ye_s \rangle. 
\end{align*}
Therefore
\begin{multline}
\label{equation: characterization - start of moment estimates}
\E*[\Q{}] \sup_{t \leq T} \left| \zeta^{\epsilon}_t(\varphi) - \int_0^t \zeta^\epsilon_s ( \mG^\dagger \varphi ) ds - \int_0^t \langle \zeta^\epsilon_s ( \varphi \overbar{h} + \alpha \overbar{\sigma}^* \nabla_x \varphi ),  d\Ye_s \rangle \right|^2 
\lesssim \E*[\Q{}] \sup_{t \leq T} \left| \mathcal{O}_\chi(\epsilon) \right|^2 + \E*[\Q{}] \sup_{t \leq T} \left| \mathcal{O}_{\wt{\chi}}(\epsilon) \right|^2 \\
+\E*[\Q{}] \sup_{t \leq T} \left|  \int_0^t \langle \rho^{\epsilon}_s ( \varphi ( h - \overbar{h} ) ), d\Ye_s \rangle \right|^2 
+\E*[\Q{}] \sup_{t \leq T} \left|  \int_0^t \langle \rho^{\epsilon}_s ( \alpha ( \sigma - \overbar{\sigma} )^* \nabla_x \varphi ), d\Ye_s \rangle \right|^2. 
\end{multline}

From the boundedness of $b, \sigma, a, h$, the growth conditions on $D^k_x \chi$ for $|k| \leq 2$, and the finite moments of all orders for $\Q{}_{(\Xe_0, \Ze_0)}$, we have that 
\[
\lim_{\epsilon \rightarrow 0} \E*[\Q{}] \sup_{t \leq T} \left| \mathcal{O}_\chi(\epsilon) \right|^2 = 0.
\]
Similarly, from the boundedness of $b, \sigma, a, h$, the growth conditions on $D^k_x \wt{\chi}$ for $|k| \leq 2$ and $b_\mathrm{I}$, and the finite moments of all orders for $\Q{}_{(\Xe_0, \Ze_0)}$, we have that 
\[
\lim_{\epsilon \rightarrow 0} \E*[\Q{}] \sup_{t \leq T} \left| \mathcal{O}_{\wt{\chi}}(\epsilon) \right|^2 = 0.
\]
Our focus now shifts to showing that 
\begin{align}
\label{equation: characterization - start of estimates for centered terms}
\lim_{\epsilon \rightarrow 0} \E*[\Q{}] \sup_{t \leq T} \left| \int_0^t \langle \rho^{\epsilon}_s ( \varphi ( h - \overbar{h} ) ), d\Ye_s \rangle \right|^2 = 0, 
\quad \text{and} \quad
\lim_{\epsilon \rightarrow 0} \E*[\Q{}] \sup_{t \leq T} \left|  \int_0^t \langle \rho^{\epsilon}_s ( \alpha ( \sigma - \overbar{\sigma} )^* \nabla_x \varphi ), d\Ye_s \rangle \right|^2 = 0.
\end{align}
Let $\psi_h(x, z) \equiv \varphi ( h - \overbar{h} )(x, z)$ and $\psi_\sigma(x, z) \equiv \alpha ( \sigma - \overbar{\sigma} )^* \nabla_x \varphi (x, z)$. 
Because $\psi_h$ and $\psi_\sigma$ are both centered with respect to $\mu_\infty(x)$ for each $x$, globally Lipschitz in $(x, z)$, and bounded in $(x, z)$, 
we let $\psi$ represent either $\psi_h$ or $\psi_\sigma$ and perform the same analysis for both. 

Applying BDG to either term in Eq. \ref{equation: characterization - start of estimates for centered terms}, we get
\begin{align}
\label{equation: characterizing limit, bound on unnormalized measure with centered term}
\E*[\Q{}] \sup_{t \leq T} \left|  \int_0^t \langle \rho^{\epsilon}_s ( \psi ), d\Ye_s \rangle \right|^2
\lesssim \E*[\Q{}] \int_0^T \left| \rho^{\epsilon}_s ( \psi ) \right|^2 ds 
= \E*[\Q{}] \int_0^T \left| \E[\Pe]{ \psi (\Xe_s, \Ze_s) \wDe_s }[\mYe_s] \right|^2 ds.
\end{align}
We now follow the argument by Kushner \cite[Chapter 6]{Kushner:1990} to partition the domain of the time integral into intervals of length at most $0 < \delta \ll 1$, where $\delta = \delta(\epsilon)$ will later be chosen as a function of $\epsilon$. 
Let $N = \left\lfloor \frac{T}{\delta} \right\rfloor \in \N{}_0$ such that $T = N \delta + \mathcal{O}(\delta)$. 
Then we have, 
\begin{align}
\E*[\Q{}] \int_0^T \left| \E[\Pe]{ \psi (\Xe_s, \Ze_s) \wDe_s }[\mYe_s] \right|^2 ds
&= \E*[\Q{}] \sum_{i = 0}^{N - 1} \int_{t_i}^{t_{i + 1}} \left| \E[\Pe]{ \psi (\Xe_s, \Ze_s) \wDe_s }[\mYe_s] \right|^2 ds \\
&+ \E*[\Q{}] \int_{N\delta}^T \left| \E[\Pe]{ \psi (\Xe_s, \Ze_s) \wDe_s }[\mYe_s] \right|^2 ds. \nonumber 
\end{align}
with $t_{i + 1} - t_i = \delta$, $\forall i$. 
We now consider a single time integral from $[t_i, t_{i+1}]$. 
For simplicity and clarity, let us use the notation $[t, t + \delta]$ instead. 
The analysis for the remainder term, over the interval $[N\delta, T]$, will follow from the same arguments. 

We introduce terms to the conditional expectation with arguments $\Xe_t$ and $\wDe_t$ fixed at the initial time of the integral over $[t, t + \delta]$, to get, 
\begin{align}
\label{equation: characterizing limit, expansion for bound on unnormalized measure with centered term}
\left| \E[\Pe]{ \psi (\Xe_s, \Ze_s) \wDe_s }[\mYe_s] \right|^2
&\lesssim \left| \E[\Pe]{ ( \psi (\Xe_s, \Ze_s) - \psi (\Xe_t, \Ze_s) ) \wDe_s }[\mYe_s] \right|^2 \\
&+ \left| \E[\Pe]{ \psi (\Xe_t, \Ze_s) ( \wDe_s - \wDe_t ) }[\mYe_s] \right|^2 
+ \left| \E[\Pe]{ \psi (\Xe_t, \Ze_s) \wDe_t }[\mYe_s] \right|^2. \nonumber
\end{align}
The first term on the right-hand side of Eq. \ref{equation: characterizing limit, expansion for bound on unnormalized measure with centered term}, by way of Jensen's inequality and Fubini on the first line, change of measure, Cauchy-Schwarz, Jensen's inequality, and the tower property of conditional expectation on the second line, contributes 
\begin{multline*}
\E*[\Q{}] \int_t^{t + \delta} \left| \E[\Pe]{ ( \psi (\Xe_s, \Ze_s) - \psi (\Xe_t, \Ze_s) ) \wDe_s }[\mYe_s] \right|^2 ds
\leq \int_t^{t + \delta} \E*[\Q{}]{} \E[\Pe]{ \left| \left( \psi (\Xe_s, \Ze_s) - \psi (\Xe_t, \Ze_s) \right) \wDe_s \right|^2 }[\mYe_s] ds \\
\leq \int_t^{t + \delta} \E[\Pe]{\left( \wDe_s \right)^2}^{1/2} \E[\Pe]{ \left| \left( \psi (\Xe_s, \Ze_s) - \psi (\Xe_t, \Ze_s) \right) \wDe_s \right|^4 }^{1/2} ds.
\end{multline*}
By Lemma \ref{lemma: moments of Girsanov}, $\E[\Pe]{\left( \wDe_s \right)^2}^{1/2} < \infty$, and by application of Cauchy-Schwarz and then the Lipschitz property of $\psi$,  we get
\begin{align*}
\E[\Pe]{ \left| \left( \psi (\Xe_s, \Ze_s) - \psi (\Xe_t, \Ze_s) \right) \wDe_s \right|^4 }^{1/2}
&\leq \E[\Pe]{ \left| \psi (\Xe_s, \Ze_s) - \psi (\Xe_t, \Ze_s) \right|^8 }^{1/4} \E[\Pe]{ \left( \wDe_s \right)^8 }^{1/4} \\
&\lesssim \E[\Pe]{  \left| \Xe_s - \Xe_t \right|^8 }^{1/4}.
\end{align*}
$\psi$ is globally Lipschitz in $x$, since each of the components of $\psi$ are either globally Lipschitz in $x$ or have a bounded derivative in $x$. 
Lemma \ref{lemma: moments on increments of slow process} gives 
\[
\E[\Pe]{ \left| \Xe_s - \Xe_t \right|^8 }^{1/4} \lesssim \left( \epsilon^8 ( 1 + \E*[\Pe]{ | \Ze_0 |^q } ) + \delta^4 ( 1 + \epsilon^8 ) ( 1 + \E{ |\Ze_0|^q } ) \right)^{1/4},
\]
for some $q \geq 0$, and therefore by the finite moments of $\Q{}_{(\Xe_0, \Ze_0)}$, the first term of Eq. \ref{equation: characterizing limit, expansion for bound on unnormalized measure with centered term} is bounded by, 
\begin{align}
\label{equation: characterization, first limit}
\E*[\Q{}] \int_t^{t + \delta} \left| \E[\Pe]{ ( \psi (\Xe_s, \Ze_s) - \psi (\Xe_t, \Ze_s) ) \wDe_s }[\mYe_s] \right|^2 ds
\lesssim \delta ( \epsilon^8 + \delta^4 ( 1 + \epsilon^8 ) )^{1/4}.
\end{align}

The second term of Eq. \ref{equation: characterizing limit, expansion for bound on unnormalized measure with centered term} similarly contributes, 
\begin{align*}
\E*[\Q{}]{} \int_t^{t + \delta} \left| \E[\Pe]{ \psi (\Xe_t, \Ze_s) ( \wDe_s - \wDe_t ) }[\mYe_s] \right|^2 ds
&\leq \int_t^{t + \delta} \E*[\Q{}]{} \E[\Pe]{ \left| \psi (\Xe_t, \Ze_s) ( \wDe_s - \wDe_t ) \right|^2 }[\mYe_s] ds \\
&\leq \int_t^{t + \delta} \E[\Pe]{\left( \wDe_s \right)^2}^{1/2} \E[\Pe]{ \left| \psi (\Xe_t, \Ze_s) ( \wDe_s - \wDe_t ) \right|^4 }^{1/2} ds.
\end{align*}
We now use the boundedness of $\psi$ to get
\begin{align*}
\E[\Pe]{ \left| \psi (\Xe_t, \Ze_s) ( \wDe_s - \wDe_t ) \right|^4 }^{1/2}
\leq | \psi |_\infty^2 \E[\Pe]{ \left| \wDe_s - \wDe_t \right|^4 }^{1/2} .
\end{align*}
Lemma \ref{lemma: moments of Girsanov} gives $\E[\Pe]{ \left| \wDe_s - \wDe_t \right|^4 }^{1/2} \lesssim \delta $ and therefore 
\begin{align}
\label{equation: characterization, second limit}
\E*[\Q{}]{} \int_t^{t + \delta} \left| \E[\Pe]{ \psi (\Xe_t, \Ze_s) ( \wDe_s - \wDe_t ) }[\mYe_s] \right|^2 ds
\lesssim \delta^2 . 
\end{align}

Recall the last term in Eq. \ref{equation: characterizing limit, expansion for bound on unnormalized measure with centered term}, 
\begin{align*}
\E*[\Q{}] \int_t^{t + \delta} \left| \E[\Pe]{ \psi (\Xe_t, \Ze_s) \wDe_t }[\mYe_s] \right|^2 ds .
\end{align*}
We first consider adding and subtracting the following term within the conditional expectation, 
\begin{align*}
\psi(\Xe_t, \widehat{Z}^{\epsilon, \Xe_t}_s) \wDe_t,
\end{align*}
where $\widehat{Z}^{\epsilon, \Xe_t}$ is the process satisfying Eq. \ref{equation: SDE, fast process with fixed slow state}, but with fixed random initial condition $x = \Xe_t$. 
Then we have
\begin{multline}
\label{equation: characterizing limit, expansion of difficult term}
\E*[\Q{}] \int_t^{t + \delta} \left| \E[\Pe]{ \psi (\Xe_t, \Ze_s) \wDe_t }[\mYe_s] \right|^2 ds 
\lesssim \E*[\Q{}] \int_t^{t + \delta} \left| \E[\Pe]{ \psi (\Xe_t, \widehat{Z}^{\epsilon, \Xe_t}_s) \wDe_t }[\mYe_s] \right|^2 ds  \\
+ \E*[\Q{}] \int_t^{t + \delta} \left| \E[\Pe]{ \left( \psi (\Xe_t, Z^{\epsilon}_s) - \psi (\Xe_t, \widehat{Z}^{\epsilon, \Xe_t}_s) \right) \wDe_t }[\mYe_s] \right|^2 ds. 
\end{multline}

Concentrating on the second term of Eq. \ref{equation: characterizing limit, expansion of difficult term},
\begin{multline*}
\E*[\Q{}] \int_t^{t + \delta} \left| \E[\Pe]{ \left( \psi (\Xe_t, Z^{\epsilon}_s) - \psi (\Xe_t, \widehat{Z}^{\epsilon, \Xe_t}_s) \right) \wDe_t }[\mYe_s] \right|^2 ds 
\lesssim \int_t^{t + \delta} \E[\Pe]{ \left| \psi (\Xe_t, Z^{\epsilon}_s) - \psi (\Xe_t, \widehat{Z}^{\epsilon, \Xe_t}_s) \right|^8 }^{1/4} ds  \\
\lesssim \int_t^{t + \delta} \E[\Pe]{ \E[\Pe]{ \left| \psi (\Xe_t, Z^{\epsilon}_s) - \psi (\Xe_t, \widehat{Z}^{\epsilon, \Xe_t}_s) \right|^8 }[\mF^{\Xe}_t \vee \mF^{\Ze}_t] }^{1/4} ds \\
\lesssim \int_t^{t + \delta} \E[\Pe]{ \E[\Pe]{ \left| \psi (x, Z^{\epsilon; (t, x, z)}_s) - \psi (x, \widehat{Z}^{\epsilon, x; (t, z)}_s) \right|^8 }[(x, z) = (\Xe_t, \Ze_t)] }^{1/4} ds . 
\end{multline*}
From the global Lipschitz property of $\psi$ in the $z$ component, we have the following estimate, 
\begin{align*}
\E[\Pe]{\left| \psi (x, Z^{\epsilon; (t, x, z)}_s) - \psi (x, \widehat{Z}^{\epsilon, x; (t, z)}_s) \right|^8}
&\lesssim \E[\Pe]{ \left| Z^{\epsilon; (t, x, z)}_s - \widehat{Z}^{\epsilon, x; (t, z)}_s \right|^8 }.
\end{align*}
In what follows, we use the notation $(X^{\epsilon; (t, x)}, Z^{\epsilon; (t, z)})$ for the pair process realized by $Z^{\epsilon; (t, x, z)}$. 
Similarly, we use $X^{\epsilon; (t, x, z)}$ when we must make clear that we are referring to the first entry of the pair $(X^{\epsilon; (t, x)}, Z^{\epsilon; (t, z)})$ which satisfies Eq. \ref{equation: multiscale signal process}.
The previous inequality is then bounded as follows, 
\begin{multline*}
\E[\Pe]{ \left| Z^{\epsilon; (t, x, z)}_s - \widehat{Z}^{\epsilon, x; (t, z)}_s \right|^8 } \lesssim \\
\frac{\delta^7}{\epsilon^{16}} \int_t^{t + \delta} \E*[\Pe] \left| f(X^{\epsilon; (t, x)}_s, Z^{\epsilon; (t, z)}_s) - f(x, \widehat{Z}^{\epsilon, x; (t, z)}_s) \right|^8 ds 
+ \frac{\delta^3}{\epsilon^8} \int_t^{t + \delta} \E*[\Pe] \left| g(X^{\epsilon; (t, x)}_s, Z^{\epsilon; (t, z)}_s) - g(x, \widehat{Z}^{\epsilon, x; (t, z)}_s)  \right|^8 ds \\
\lesssim \frac{\delta^7}{\epsilon^{16}} \int_t^{t + \delta} \E*[\Pe] \left| f(X^{\epsilon; (t, x)}_s, Z^{\epsilon; (t, z)}_s) - f(x, Z^{\epsilon; (t, z)}_s) \right|^8 
+ \E*[\Pe] \left| f(x, Z^{\epsilon; (t, x, z)}_s) - f(x, \widehat{Z}^{\epsilon, x; (t, z)}_s) \right|^8 ds \\
+ \frac{\delta^3}{\epsilon^8} \int_t^{t + \delta} \E*[\Pe] \left| g(X^{\epsilon; (t, x)}_s, Z^{\epsilon; (t, z)}_s) - g(x, Z^{\epsilon; (t, z)}_s) \right|^8 
+ \E*[\Pe] \left| g(x, Z^{\epsilon; (t, x, z)}_s) - g(x, \widehat{Z}^{\epsilon, x; (t, z)}_s)  \right|^8 ds \\
\leq \frac{\delta^3}{\epsilon^8} \left( \frac{\delta^4}{\epsilon^8} | \nabla_x f |_\infty^8 + | \nabla_x g |_\infty^8 \right) \int_t^{t + \delta} \E*[\Pe] \left| X^{\epsilon; (t, x, z)}_s - x \right|^8 ds \\
+ \frac{\delta^3}{\epsilon^8} \left( \frac{\delta^4}{\epsilon^8} | \nabla_z f |_\infty^8 + | \nabla_z g |_\infty^8 \right) \int_t^{t + \delta} \E*[\Pe] \left| Z^{\epsilon; (t, x, z)}_s - \widehat{Z}^{\epsilon, x; (t, z)}_s \right|^8 ds.
\end{multline*}
From Lemma \ref{lemma: moments on increments of slow process}, for some $q \geq 0$, we get 
\[
\int_t^{t + \delta} \E*[\Pe] \left| X^{\epsilon; (t, x, z)}_s - x \right|^8 ds \lesssim \delta \epsilon^8 ( 1 + |z|^q ) + \delta^5 ( 1 + \epsilon^8 ) ( 1 + |z|^q ).
\] 
Let 
\[
\eta(\epsilon, \delta) \equiv \left( \frac{\delta^8}{\epsilon^{16}} + \frac{\delta^4}{\epsilon^8} \right). 
\]
Therefore Gr\"onwall gives us
\begin{align*}
\E[\Pe]{ \left| Z^{\epsilon; (t, x, z)}_s - \widehat{Z}^{\epsilon, x; (t, z)}_s \right|^8 }
\lesssim \eta(\epsilon, \delta) \left( \epsilon^8 + \delta^4 ( 1 + \epsilon^8 ) \right) \exp( \eta(\epsilon, \delta) ) ( 1 + |z|^q ).
\end{align*}
For further brevity, let us define
\[
\frak{F}(\epsilon, \delta) 
\equiv 
\eta(\epsilon, \delta) \left( \epsilon^8 + \delta^4 ( 1 + \epsilon^8 ) \right) \exp( \eta(\epsilon, \delta) ).
\]
Therefore the second term in Eq. \ref{equation: characterizing limit, expansion of difficult term} is bounded by
\begin{align}
\label{equation: characterization, third limit}
\E*[\Q{}] \int_t^{t + \delta} \left| \E[\Pe]{ \left( \psi (\Xe_t, Z^{\epsilon}_s) - \psi (\Xe_t, \widehat{Z}^{\epsilon, \Xe_t}_s) \right) \wDe_t }[\mYe_s] \right|^2 ds 
\lesssim \int_t^{t + \delta} \E[\Pe]{ \frak{F}(\epsilon, \delta) ( 1 + |\Ze_t|^q ) }^{1/4} ds \nonumber \\ 
\lesssim \delta \frak{F}(\epsilon, \delta)^{1/4} ( 1 + \E[\Pe]{ |\Ze_0|^{q'} } )^{1/4} 
\lesssim \delta \frak{F}(\epsilon, \delta)^{1/4} .
\end{align}
For the first term on the right hand side of Eq. \ref{equation: characterizing limit, expansion of difficult term}, we condition the centering term on a larger filtration $\mathcal{H} = \mYe_s \vee \mF^{\Xe}_t \vee \mF^{\Ze}_t$, and then use the fact that $\sigma(\widehat{Z}^{\epsilon, \Xe_t}_s) \vee \mYe_t \vee \mF^{\Xe}_t \vee \mF^{\Ze}_t$ is independent of $\sigma(\Ye_r - \Ye_t; r \in [t, s])$ under $\Pe$ and that $(\Xe_t, \widehat{Z}^{\epsilon, \Xe_t}_s)$ is Markov in the larger filtration $\mYe_s \vee \mF^{\Xe}_t \vee \mF^{\Ze}_s$ to yield,
\begin{align*}
\E*[\Q{}] \int_t^{t + \delta} \left| \E[\Pe]{ \psi (\Xe_t, \widehat{Z}^{\epsilon, \Xe_t}_s) \wDe_t }[\mYe_s] \right|^2 ds
&= \E*[\Q{}] \int_t^{t + \delta} \left| \E[\Pe]{ \E[\Pe]{ \psi (\Xe_t, \widehat{Z}^{\epsilon, \Xe_t}_s) }[\mathcal{H}] \wDe_t }[\mYe_s] \right|^2 ds \\
&= \E*[\Q{}] \int_t^{t + \delta} \left| \E[\Pe]{ \E[\Pe]{ \psi (x, \widehat{Z}^{\epsilon, x; (t, z)}_s) }[(x, z) = (\Xe_t, \Ze_t)] \wDe_t }[\mYe_s] \right|^2 ds.
\end{align*}
Applications of Jensen's inequality, Cauchy-Schwarz, the tower property, Lemma \ref{lemma: moments of Girsanov} and \ref{lemma: bounds on centered functions under semigroup} then give the estimate, 
\begin{multline}
\label{equation: characterization, fourth limit}
\E*[\Q{}] \int_t^{t + \delta} \left| \E[\Pe]{ \E[\Pe]{ \psi (x, \widehat{Z}^{\epsilon, x; (t, z)}_s) }[(x, z) = (\Xe_t, \Ze_t)] \wDe_t }[\mYe_s] \right|^2 ds \\
\lesssim \int_t^{t + \delta} \E[\Pe]{ \left| \E[\Pe]{\psi (x, \widehat{Z}^{\epsilon, x; (t, z)}_s)}[(x, z) = (\Xe_t, \Ze_t)] \right|^8 }^{1/4} ds 
= \int_t^{t + \delta} \E[\Pe]{ \left| T^{F, x}_{(s - t) / \epsilon^2}(\psi(\Xe_t, \cdot))(\Ze_t ) \right|^8 }^{1/4} ds  \\
= \epsilon^2 \int_0^{\delta / \epsilon^2} \E[\Pe]{ \left| T^{F, x}_u(\psi(\Xe_t, \cdot))(\Ze_t ) \right|^8 }^{1/4} du 
\leq \epsilon^2 \int_0^\infty \E[\Pe]{ \left| T^{F, x}_u(\psi(\Xe_t, \cdot))(\Ze_t ) \right|^8 }^{1/4} du \\
\lesssim \epsilon^2 \left( 1 + \E[\Pe]{ \left| \Ze_t \right|^q } \right)^{1 / 4} 
\lesssim \epsilon^2 \left( 1 + \E[\Q{}]{ \left| \Ze_0 \right|^{q'} } \right)^{1 / 4} 
\lesssim \epsilon^2 .
\end{multline}

Collecting all our bounds for Eq. \ref{equation: characterizing limit, bound on unnormalized measure with centered term}, that is Eqs. \ref{equation: characterization, first limit}, \ref{equation: characterization, second limit}, \ref{equation: characterization, third limit}, and \ref{equation: characterization, fourth limit}, and accounting for the discretization of the time integral into $N$ segments, which results in $T / \delta$ times the estimates, we have 
\begin{align}
\label{equation: characterizing limit, collection of all terms - sensor function}
\E*[\Q{}] \sup_{t \leq T} \left|  \int_0^t \langle \rho^{\epsilon}_s ( \psi ), d\Ye_s \rangle \right|^2
\lesssim ( \epsilon^8 + \delta^4 ( 1 + \epsilon^8 ) )^{1/4}
+ \delta 
+ \frak{F}(\epsilon, \delta)^{1/4}  
+ \frac{\epsilon^2}{\delta} . 
\end{align}
If we choose $\delta(\epsilon) = \epsilon^2 (- \ln \epsilon )^p$ with $p \in (0, 1/8)$, then $\lim_{\epsilon \rightarrow 0^{+}} \delta(\epsilon) = 0$ and $\frak{F}(\epsilon, \delta) \rightarrow 0$ (see Lemma \ref{proposition: limit result for characterization proof (qualitative convergence)}), which completes the proof. 
\end{proof}
}
}

{
\lem{
\label{lemma: uniqueness of weak limits}
Under either of the assumptions: 
\begin{enumerate}[a.]
\item the coefficients of $\mG^\dagger$ and $\overbar{h}, \overbar{\sigma}$ are $C^{2 + \alpha}_b$, for some $\alpha \in (0, 1)$, or
\item $\overbar{a} + \wt{a} \succ 0$ uniformly in $x$ and the coefficients of $\mG^\dagger$ and $\overbar{h}, \overbar{\sigma}$ are $C^\alpha_b$, for some $\alpha \in (0, 1)$,
\end{enumerate}
the finite signed Borel measure-valued process $\zeta$, has the unique solution $\zeta_t = 0$, $\Q{}$-a.s. $\forall t \in [0, T]$.
\begin{proof} 
Our objective is simply to show that 
\begin{align*}
\zeta_t(\varphi) = \int_0^t \zeta_s ( \mG^\dagger \varphi ) ds + \int_0^t \langle \zeta_s ( \varphi \overbar{h} + \alpha \overbar{\sigma}^* \nabla_x \varphi ), dY_s \rangle, \quad \zeta_0(\varphi) = 0,
\end{align*}
is a Zakai equation, since uniqueness then follows from \cite[Theorem 3.1, p.454]{Rozovskii:1991}. 

Let $X^0$ be the diffusion process with infinitesimal generator $\mG^\dagger$. 
In particular, consider the following system of equations, 
\begin{align}
\tag{\ref{equation: SDE to achieve averaged unnormalized filter equation}}
dX^0_t 
&= \left[ \overbar{b}(X^0_t) + \wt{b}(X^0_t) \right] dt + \wt{a}^{1/2}(X^0_t) d\wt{W}_t + ( \overbar{a}(X^0_t) - \overbar{\sigma} \overbar{\sigma}^*(X^0_t) )^{1/2} d\widehat{W}_t + \overbar{\sigma}(X^0_t) dW_t \\
dY_t 
&= 
\overbar{h}(X^0_t) dt + \alpha dW_t + \gamma dB_t,  \nonumber
\end{align}
where $\alpha dW_t + \gamma dB_t$ is a standard Brownian motion, $\wt{W}, \widehat{W}, W, B$ are independent standard Brownian motions under $\Q{}$. 
This system of equations yield a Zakai equation of the desired form after the change of measure given by $D_t = \exp( - \int_0^t \langle \overbar{h}(X^0_s), \alpha dW_s + \gamma dB_s \rangle - \frac{1}{2} \int_0^t | \overbar{h}(X^0_s) |^2 ds)$ is performed. 
\FutureWork{
As Nicolas points out, application of Jensen's inequality shows that indeed the term in question is semi-positive definite, and therefore we do not have an issue taking the square root. 
\begin{multline*}
|z \cdot E[\sigma] E[\sigma*] z|
= |\sum_k E[ \sum_i z_i \sigma_{i k}] E[ \sum_j z_j \sigma_{j k}]|
=  \le \sum_k E[ |\sum_i z_i \sigma_{i k}|^2 ] \\
= \sum_{i,j} z_i \sum_k E[\sigma_{ik} \sigma_{jk}] z_j
= \sum_{i,j} z_i E[a_{ij}] z_j
= z \cdot E[a] z.
\end{multline*}
A separate issue though is whether the square root is smooth. 
For this we need to look at Stoock Lemma 2.3.3. 
Nicolas states that if the term in question is $C^2$ with uniformly bounded second derivative, then we will be OK. 
Additional remarks that must be made, $V$ may be of smaller dimension than $U$.
Need to make the observation in Section 2 that $\overbar{\sigma}$ may be equal to zero and therefore the SDE we get in the end may have less correlation.  
}

\end{proof}
}
}

{
\thm{
\label{theorem: weak convergence of unnormalized filter}
Assume that $f$ and $g$ satisfy \ref{assumption: positive recurrence} and \ref{assumption: uniform ellipticity}, that $b_\mathrm{I}$ is centered with respect to $\mu_\infty(x)$ for each $x$ and that $\Q{}_{(\Xe_0, \Ze_0)}$ has finite moments of every order.
Additionally, assume either:
\begin{enumerate}[a.]
\item \hyperlink{H}{$H^{3, 2 + \alpha}$} holds for $\alpha \in (0, 1)$;
for each $z$, $b(\cdot, z) , \sigma(\cdot, z) \in C^3$, and $b_\mathrm{I}(\cdot, z) \in C^4$;
that $b$ and $b_\mathrm{I}$ are Lipschitz in $z$, and $\sigma$ is globally Lipschitz in $z$; 
that $b, b_\mathrm{I}, \sigma$ satisfy the growth conditions
\begin{align*}		
\left| b(x, z) \right| + \left| b_\mathrm{I}(x, z) \right| + \left| \sigma \sigma^* (x, z) \right| &\leq C (1 + |z|)^\beta, \\
\sum_{|k| = 1}^2 \left| D^k_x b(x, z) \right| +\left| D^k_x \sigma \sigma^* (x, z) \right| &\leq C (1 + |z|^q), \\
\sum_{|k| = 1}^3 \left| D^k_x b_\mathrm{I}(x, z) \right| &\leq C (1 + |z|^q),
\end{align*}
for some \Nicolas{$\beta < -2$} and $q > 0$;
that $h$ is bounded in $(x, z)$, $h(\cdot, z) \in C^3$ for each $z$, and $h$ is globally Lipschitz in $z$.
\item $\overbar{a} + \wt{a} \succ 0$ uniformly in $x$;
\hyperlink{H}{$H^{2, 2 + \alpha}$} holds for $\alpha \in (0, 1)$;
for each $z$, $b(\cdot, z), b_\mathrm{I}(\cdot, z), \sigma(\cdot, z) \in C^2$; 
that $b$ and $b_\mathrm{I}$ are Lipschitz in $z$, and $\sigma$ is globally Lipschitz in $z$;
that $b, b_\mathrm{I}, \sigma$ satisfy the growth conditions
\begin{align*}		
\left| b(x, z) \right| + \left| b_\mathrm{I}(x, z) \right| + \left| \sigma \sigma^* (x, z) \right| &\leq C (1 + |z|)^\beta, \\
\sum_{|k| = 1}^2 \left| D^k_x b(x, z) \right| + \left| D^k_x b_\mathrm{I}(x, z) \right| + \left| D^k_x \sigma \sigma^* (x, z) \right| &\leq C (1 + |z|^q),
\end{align*}
for some \Nicolas{$\beta < -2$} and $q > 0$;
$h$ is bounded in $(x, z)$, that $h$ is globally Lipschitz in $(x, z)$. 
If $a \succ 0$, which implies $\overbar{a} + \wt{a} \succ 0$, then the Lipschitz condition in $z$ for $b, b_\mathrm{I}$ can be relaxed to $\alpha$-H\"older continuity. 
\end{enumerate}
Then $\zeta^\epsilon = \rho^{\epsilon, x} - \rho^0 \Rightarrow 0$ as $\epsilon \rightarrow 0$. 
\begin{proof}
This follows from Corollary \ref{corollary: solutions to signed Zakai equations are relatively compact}\textendash the existence of weak limits of the probability measures induced on path space by $\zeta^\epsilon$, Lemma \ref{lemma: characterization of weak limits}\textendash the characterization of the limit points, and Lemma \ref{lemma: uniqueness of weak limits} on the uniqueness of the limiting evolution equation.
\end{proof}
}
}

{
\lem{
\label{lemma: weak convergence of unnormalized filter implies convergence of normalized}
Let $\rho^\epsilon$ be a solution of Eq. \ref{equation: Zakai} and $\rho^0$ a solution of Eq. \ref{equation: homogenized Zakai}.
Assume that $h, \overbar{h}$ and the coefficients of $\mG^\dagger$ are bounded. 
If $\rho^{\epsilon, x} - \rho^0 \Rightarrow 0$ as $\epsilon \rightarrow 0$, then $\pi^{\epsilon, x} - \pi^0 \Rightarrow 0$. 
\begin{proof} 
Let $\varphi \in C^2_b(\R{m}; \R{})$ and $t \in [0, T]$, then 
\begin{align*}
(\pi^{\epsilon, x} - \pi^0)_t (\varphi)  
= \frac{\rho^{\epsilon, x}_t(\varphi)}{\rho^{\epsilon, x}_t(1)} - \frac{\rho^0_t(\varphi)}{\rho^0_t(1)} 
&= \frac{ (\rho^{\epsilon, x} - \rho^0)_t(\varphi) }{ \rho^{\epsilon, x}_t(1) } + \pi^0_t(\varphi) \frac{ ( \rho^0 - \rho^{\epsilon, x})_t(1) }{ \rho^{\epsilon, x}_t(1) }.
\end{align*}
The weak convergence of $(\pi^{\epsilon, x} - \pi^0)_t$ now follows from the estimate
\begin{align*}
\lim_{\delta \rightarrow 0} \inf_{\epsilon > 0} \Q{} \left( \inf_{t \leq T} \rho^{\epsilon, x}_t(1) > \delta \right) = 1,
\end{align*}
and the fact that $\varphi$ is bounded and $\pi^0_t$ is almost surely equal to a probability measure.
\end{proof}
}
}

\section{Remark on Conditions for the Fast Semigroup}
\label{section: remark on condition for the fast semigroup}

The necessary conditions in this paper are sometimes at odds with Theorems 2 and 3 from \cite[p.1171]{Pardoux:2003}, which are used in this paper for a number of propositions and theorems listed below. 
Specifically, in \cite{Pardoux:2003}, the condition in Theorems 2 and 3 are given as \hyperlink{H}{$H^{1, 2 + \alpha}$} (there are actually two scenarios to consider, but in this paper we only consider one of them, which is the one just quoted). 
In particular, only one continuous derivative in the $x$-component is ever needed in the coefficients $f$ and $g$ to be able to take $k \geq 1$ derivatives of the new function under the semigroup $T^{F, x}(\varphi)$, where $\varphi \in C^k$ for instance. 
Because the Poisson solution of \cite[Theorem 3, p.1171]{Pardoux:2003} is proven based on Theorem 2, the same condition of \hyperlink{H}{$H^{1, 2 + \alpha}$} shows up there, even if $k \geq 1$ derivatives of the Poisson solution are desired. 
A counter example as to why this condition is insufficient is given next.

\subsection{Counter Example}

Let $g(x, z) = g(x)$ depend only on $x$ and let $f(x, z) = -z$. 
Then the fast process is 
\[
dZ^x_t = -Z^x_t dt + g(x) dB_t, 
\]
and therefore $Z^x$ is an Ornstein-Uhlenbeck process (in particular Gaussian), and hence satisfies the recurrence condition for \cite[p.1171]{Pardoux:2003}. 
We can choose $g$ to satisfy the uniform ellipticity condition as well, assume this to be true. 
If $Z^x_0 = z$, then 
\[ 
Z^x_t = e^{- t} z + \int_0^t e^{- (t - s)} g (x) d B_s \sim \mathcal{N} \left( e^{- t} z, \frac{g (x)^2}{2} (1 - e^{- 2 t}) \right), 
\]
and thus the transition density at time $t$ in $z$, having started from $(x, z')$ at the initial time is 
\[ 
p_t (z, z' ; x) = r \left( \frac{g (x)^2}{2} (1 - e^{- 2 t}), z' -  e^{- t} z \right), 
\]
where $r (s, y)$ is the Gaussian density with variance $s$, evaluated in $y$.
Consider now the test function $\psi (x, z) = \cos (z)$, which is infinitely smooth in $x$ (and in $z$). 
Note that for $Y \sim \mathcal{N} (\mu, g^2)$ we have
\[ 
\mathbb{E} [\cos (Y)] = \frac{1}{2} \mathbb{E} [e^{i Y} + e^{- i Y}] 
= \frac{1}{2} \left( e^{i \mu - \frac{1}{2} g^2} + e^{- i \mu -\frac{1}{2} g^2} \right) 
= e^{- \frac{1}{2} g^2} \cos (\mu), 
\]
and therefore the semigroup (notation from \cite[p.1171]{Pardoux:2003}) is 
\[ 
p_t (z, \psi ; x) =\mathbb{E}_z [\cos (Z^x_t)] 
= \exp \left( - \frac{1}{2} \frac{g (x)^2}{2} (1 - e^{- 2 t}) \right) \cos (e^{- t} z) . 
\]
If $g^2 \notin C^2$, then this function is not $C^2$ in $x$.

\subsection{List of Changes}

The condition should be \hyperlink{H}{$H^{k, 2 + \alpha}$}, and we use this condition instead of the one given in \cite{Pardoux:2003}. 
The difference in the requirements of various propositions and theorems are subtle, but listed here for reference: 
\begin{enumerate}[(i)]
\item Theorem \ref{theorem: pardoux and veretennikov diffusion approximation}, \hyperlink{H}{$H^{1, 2 + \alpha}$} has become \hyperlink{H}{$H^{2, 2 + \alpha}$}.
\item Theorem \ref{theorem: solution, bounds and regularity of Poisson equation}, \hyperlink{H}{$H^{1 , 2 + \alpha}$} has become \hyperlink{H}{$H^{k , 2 + \alpha}$}.
\item Lemma \ref{lemma: regularity and integrability of poisson averaged coefficients}, \hyperlink{H}{$H^{1, 2 + \alpha}$} has become \hyperlink{H}{$H^{j, 2 + \alpha}$}.
\item Theorem \ref{theorem: weak convergence of unnormalized filter}, for a.) \hyperlink{H}{$H^{3, 2 + \alpha}$} has become \hyperlink{H}{$H^{4, 2 + \alpha}$}. This was a result of needing the third derivative in $x$ of $\wt{b}$, which required the fourth derivative in $x$ of the Poisson solution $\mG^{-1}_F(b_\mathrm{I})$. 
\end{enumerate}

A final remark, is that Lemma \ref{lemma: bounds on centered functions under semigroup} is not affected by this, because there we are also using the density result of \cite[Theorem 1, p.1170]{Pardoux:2003}, which is correct and requires stronger conditions than \cite[Theorem 2, p.1171]{Pardoux:2003}.

\section*{Acknowledgement}
R.B. and N.S.N. acknowledge partial support for this work from the Air Force Office of Scientific Research under grant number FA9550-17-1-0001, and N.S.N. acknowledges partial support from the National Sciences and Engineering Research Council Discovery grant 50503-10802.

\setcounter{section}{0}
\renewcommand{\thesection}{\Alph{section}}
\section{Appendix}

{
\lemma{
\label{proposition: limit result for characterization proof (qualitative convergence)}
Let $p \in (0, 1/8)$, $\delta(\epsilon) = \epsilon^2 (- \ln \epsilon )^p$, then 
\[
\lim_{\epsilon \rightarrow 0^{+}} \left( \frac{\delta^8}{\epsilon^{16}} + \frac{\delta^4}{\epsilon^8} \right) (\epsilon^8 + \delta^4  (1 + \epsilon^8)) \exp \left( \frac{\delta^8}{\epsilon^{16}} + \frac{\delta^4}{\epsilon^8} \right) = 0.
\]
\begin{proof}
We first expand the expression with the choice of $\delta(\epsilon)$ to get, 
\[
\left( (-\ln \epsilon)^{8p} + (-\ln \epsilon)^{4p} \right) (\epsilon^8 + \epsilon^8 (-\ln \epsilon)^{4p} + \epsilon^{16} (-\ln \epsilon)^{4p}) \exp \left( (-\ln \epsilon)^{8p} + (-\ln \epsilon)^{4p} \right).
\]
Expanding and distributing the terms, we identify the term that would be most limiting for convergence to zero, 
\begin{align*}
(-\ln \epsilon)^{12p} \epsilon^8 \exp \left( (-\ln \epsilon)^{8p} + (-\ln \epsilon)^{4p} \right)
\lesssim \epsilon^7 \exp \left( 2 (-\ln \epsilon )^{8p} \right). 
\end{align*}
Since $8p < 1$, for all sufficiently small $\epsilon > 0$ we have, 
\[
\exp \left( 2 (-\ln \epsilon )^{8p} \right) \leq \exp \left( - 2 \ln \epsilon \right) = \epsilon^{-2},
\]
and therefore 
\[
\lim_{\epsilon \rightarrow 0^+} \epsilon^7 \exp \left( 2 (-\ln \epsilon )^{8p} \right)
\leq \lim_{\epsilon \rightarrow 0^+} \epsilon^5
= 0. 
\]
\end{proof}
}
}

\printbibliography[heading=bibintoc]
\end{document}